\documentclass{article}

\usepackage{amsmath,amssymb,graphicx} 
\usepackage[margin=0.9in]{geometry} 
\usepackage{enumerate}
\usepackage[toc]{appendix}
\usepackage{hyperref}
\usepackage{tikz}
\DeclareRobustCommand\full{\tikz[baseline=-0.6ex]\draw[thick] (0,0)--(0.5,0);}

\DeclareRobustCommand\dashed{\tikz[baseline=-0.6ex]\draw[thick,dashed] (0,0)--(0.54,0);}

\DeclareRobustCommand\dashdot{\tikz[baseline=-0.6ex]\draw[thick,dash dot] (0,0)--(0.5,0);}
\usepackage{comment}
\usepackage{diagbox}
\usepackage{multirow,multicol}
\usepackage{colortbl}
\usepackage{pgfplots}
\pgfplotsset{compat=1.8}

\begin{document}

\nocite{*} 

\title{A moving boundary approach of capturing diffusants penetration into rubber: FEM approximation and comparison with laboratory measurements}
 \author{S.\,Nepal $^{1,*}$,  R.\,Meyer$^{2}$, N.\,H.\,Kr\"oger$^{2,3}$, T.\,Aiki$^4$, A.\,Muntean$^1$, Y.\,Wondmagegne$^1$, U.\,Giese$^2$ \\
$^1$ Department of Mathematics and Computer Science, Karlstad University, Sweden\\
$^2$ Deutsches Institut f\"ur Kautschuktechnologie e.V.\,(DIK), Hanover, Germany\\
$^3$ material prediction GmbH, Wardenburg, Germany\\
$^4$ Department of Mathematics, Japan Women's University, Tokyo,  Japan\\
* surendra.nepal@kau.se}

\date{\today} 

\maketitle
\noindent
\begin{abstract}

We propose a moving-boundary scenario to model the penetration of diffusants into dense and foamed rubbers. The presented modelling approach recovers experimental findings related to the diffusion of cyclohexane and the resulting swelling in a piece of material made of ethylene propylene diene monomer rubber (EPDM). The main challenge is to find out relatively simple model components which can mimic the mechanical behavior of the rubber. Such special structure is identified here so that the computed penetration depths of the diffusant concentration are within the range of experimental measurements. We investigate two cases: a dense rubber and a rubber foam, both made of the same matrix material. 

After a brief discussion of scaling arguments, we present a finite element approximation of the moving boundary problem. To overcome numerical difficulties due to the \textit{a priori} unknown motion of the diffusants penetration front, we transform the governing model equations from the physical domain with moving unknown boundary to a fixed fictitious domain. We then solve the transformed equations by the finite element method and explore the robustness of our approximations with respect to relevant model parameters. Finally, we discuss numerical estimations of the expected large-time behavior of the material. 


\end{abstract}

\section{Introduction}


Rubber is one of most usable polymers in our daily life, where it
is employed for multiple purposes ranging from tires to materials for daily use,  medical devices, and toys. The source of the natural rubber is  the latex of rubber trees (esp. {\em hevea brasiliensis}). Due to its softness, the use of natural rubber is limited. Besides natural rubber, synthetic rubbers are widely used. They are usually produced from petroleum byproducts. To develop the wanted toughness and strength, crosslinking of the rubber molecules, e.g.\,by vulcanization, is required.  There are several ways to crosslink rubber, but the most common way is to add sulfur components that help forming bonds  between rubber molecules. This leads to  an increased strength as well as  to a lower solubility in solvents; see e.g.\,\cite{morton2013rubber} and compare also \cite{sombatsompop2000effects} and references cited therein. On the other hand, due to the solubility of low molecular components, depending on their polarity, size of the molecules as well on the flexibility of the polymer chains,  the material is to some extent permeable. This is precisely where our interest lies. Predicting the penetration of small molecules in polymeric materials  has a great importance in polymer science and engineering especially if one has in mind applications in areas such as food packaging \cite{reynier2002migration}, controlled drug delivery \cite{mircioiu2019mathematical}, membrane separation \cite{van1992diffusional}, etc.  Especially, the {\em a priori}  knowledge of critical diffusants and their distribution in rubber matrices is of great interest,  cf.\cite{Giese2000,Giese1998,Rosca2004,Rosca2006}.
To determine  the approximate position of the  penetration front, many different type of models involving various degrees of complexity have been proposed earlier; see e.g. \cite{Fasano, masoodi2012numerical, neff2019modelling, Hans, Wilmers}. However, a general consensus is not yet reached in the literature. Consequently, many fundamental questions are still open particularly what concerns the large-time behavior of rubber-based materials when they are exposed to environmental conditions. For the scenario we have in view, this is due to the fact that one does not understand yet sufficiently well how the macroscopic rheology of the rubber changes particularly when diffusants-induced swelling takes place.   

In this framework, we deal with  the macroscopic modeling and numerical simulation for the penetration of diffusants via absorption and diffusion through swellable rubbers. We consider two cases: (i) dense rubber and (ii) foamed rubber; both being made based on the same polymer chemistry. 

We are mainly addressing the following questions:
\begin{itemize}
\item[(Q1)] How far does the concentration of a diffusant species\footnote{We consider a population of small-size particles able to get absorbed into the rubber's skeleton and diffuse from there inside the material. We refer here to this diffusing population as "diffusants".} penetrate into a  dense rubber?
\item[(Q2)] How far are these diffusants able to penetrate a foamed rubber?
\end{itemize}
Although these questions may appear fairly simple, they are quite hard to answer in a full generality. We consider in this paper an experimental setup (see Section \ref{lab}) that can be traced numerically by means of one dimensional evolution model with a freely evolving moving boundary - the penetration front of the diffusants (see Section \ref{model}). We address here  (Q1) and (Q2)  by comparing the model output with the measured penetration depths and discuss the expected large-time behavior of both materials when drastic variations of parameters are applied. 

The experimental setup, detailed in Section \ref{lab}, is sketched in Figure \ref{Fig:1}. The model proposed in Section \ref{model} is meant to fit to this precise setup. As a consequence of the material's shape, the modeling domain is a 1D line (located right in the middle of the rubber sheet), which is supposed to swell as time elapses, due to a continuous penetration of diffusants up to a location at time $t>0$ denoted by $s(t)$. Predicting the time evolution of this moving front $s(t)$ is our main target - the evolution takes place along the longitudinal line indicated in in Figure \ref{Fig:1} so that  $0<s(0)\leq s(t)\leq \ell$, where $\ell$ is a characteristic length of the rubber sample.
In Section \ref{model}, we  describe the governing equations of our model. They are taken over from \cite{kumazaki2020global}, where the model was designed to suit our purpose. We start off with performing the non-dimensionalization of the model equations to identify the typical sizes of the characteristic times of diffusion, absorption and swelling, and consequently also the typical size of the mass transfer Biot number. As next step, for the ease of handling the numerical simulations and parameter identification procedure, we transform  the non-dimensional model equations from the {\em a priori} unknown  domain confined by the target diffusants's penetration depth $s(t)$ into a fixed domain.
In Section \ref{FEM}, we introduce a finite element approximation of our model, while  in Section \ref{simulation}, we
present and discuss the simulation results. We calibrate our model based on the experimental data corresponding to the dense rubber case. Afterwards in Section \ref{Secfoam}, we use the same model and capture the measured diffusion fronts in the foamed rubber by changing either the kinetic parameters (entering the speed law for the diffusion front) or the effective diffusivity), or a suitable combination of these parameters.  We close the paper with a series of discussions and hints to further work, see Section \ref{discussion} and Section \ref{outlook}. Of particular interest for us is the study of the expected large-time behavior of such materials when they are exposed to less controlled environmental conditions. This work can be seen as a preliminary step in this direction.

\section{Experimental investigations}\label{lab}

The modelling idea is based on the following experimental results related to ethylene propylene diene monomer rubber (EPDM) being swollen in cyclohexane. Especially, the difference in the diffusion respectively the swelling behaviour between full matrix rubber and its foam is of interest.\\

{\bf Material preparation:}
The material's ingredients, see Table\,\ref{Tab:Material}, were admixed at a temperature of 50$^\circ$C and rotor speed of 30 rotation per minute. Each mixture was vulcanized as plates for 24 minutes at 160$^\circ$C in a press. The full matrix material was produced in 1 mm plates, whereas plates of the foam were produced releasing the pressure of the vulcanisation press resulting in thicker plates. Subsequently, the foam plates were cut in 1 mm slices. The geometries used were 1 mm thick plates with 17 mm width and 40 mm height.

\begin{table}[h]
	\centering
	\begin{tabular}{|>{\columncolor[HTML]{C0C0C0}}c|c|c|c|c|c|c|c|}
		\hline
		Ingredients in phr & EPDM & ZnO & Stearic Acid & TAIC (70\%) & DCP (40\%) & TMQ & Backing powder \\ \hline
		Full rubber matrix & 100 & 2.5 & 2.5 & 0.6 & 8 & 0 & 0 \\ \hline
		Foam & 100 & 2.5 & 2.5 & 0.6 & 8 & 1.5 & 4.7 \\ \hline
	\end{tabular}
	\caption{Recipe of the investigated EPDM mixture (phr $\widehat{=}$ parts per hundred rubber related to mass parts).}
	\label{Tab:Material}
\end{table}

\begin{figure}[h!]
	\centering
	\includegraphics[height=6.5cm]{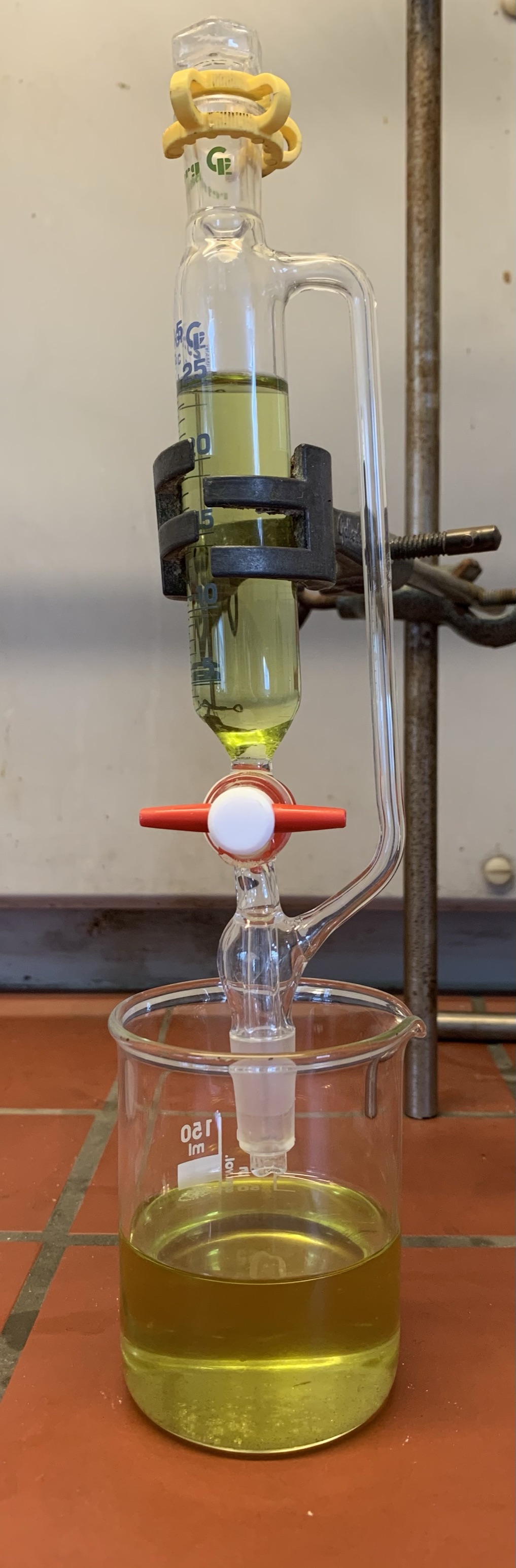}
	\hspace{0.1cm}
	\includegraphics[height=6.5cm]{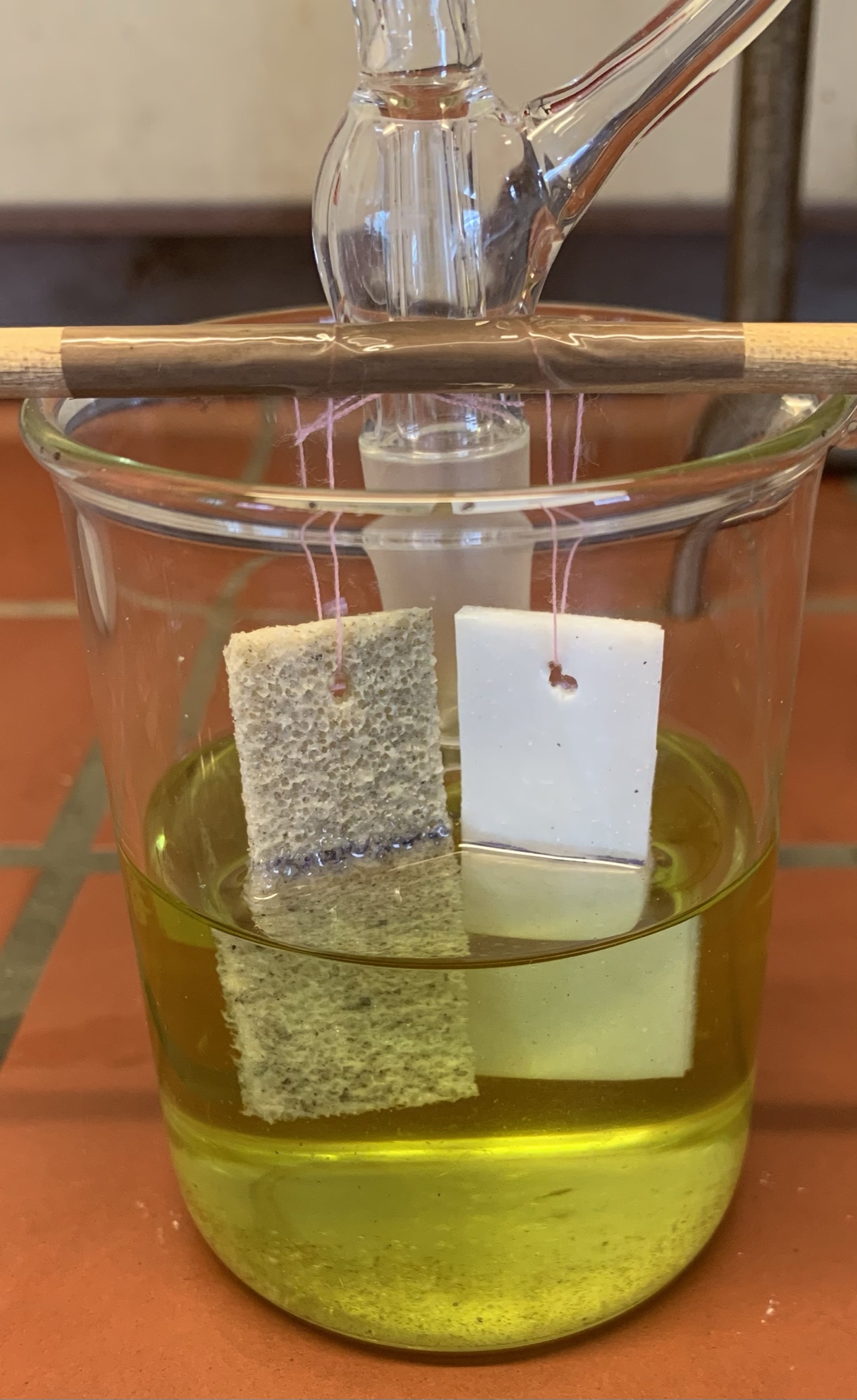}
	\hspace{0.1cm}
	\includegraphics[height=6.5cm]{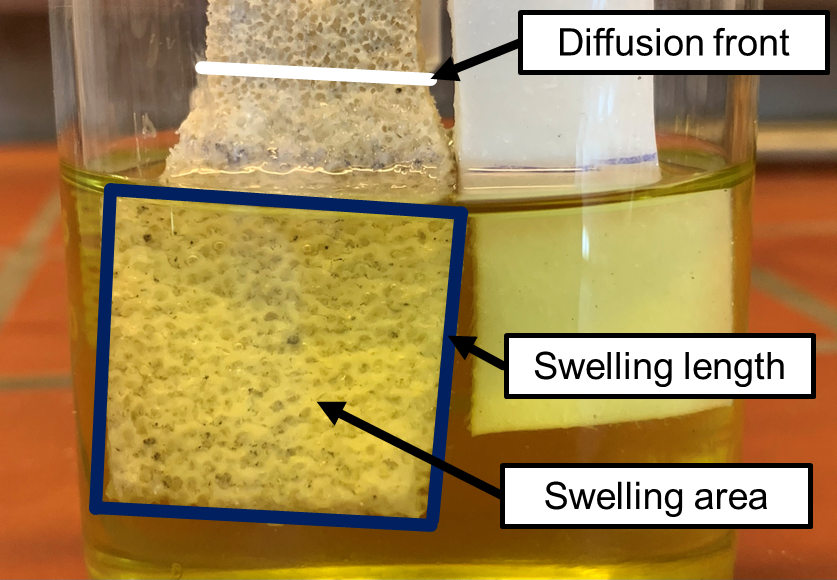}
	\caption{Experimental setup and evaluation.}
	\label{Fig:Exp}
\end{figure}

\begin{figure}[h!]
    \centering
   \begin{tikzpicture}
	\pgfplotsset{scale only axis,scaled x ticks=base 10:0,xmin=, 		xmax=330}
	\begin{axis}[
	axis y line*=left,
	ymin=0, ymax=8,
	xlabel=Time in min,
	ylabel=Diffusion front in mm,
	domain=0:350,
	grid=both,
	grid style={line width=.1pt, draw=gray!10},
	major grid style={line width=.2pt,draw=gray!500}]
	\pgfplotsset{every axis legend/.append style={
			at={(-0.05,-0.2)},
			anchor=north west}}
	\addplot [mark=*,black,very thick] table [x=t, y=FR, col sep=semicolon]  {Exp.csv};
	\addlegendentry{Rubber (front)};
	\addplot [mark=*,red,very thick] table [x=t, y=FF, col sep=semicolon]  {Exp.csv};
	\addlegendentry{Foamed rubber (front)};
	\end{axis}
	\begin{axis}[
	axis y line*=right,
	ylabel=Submerged length in mm,
	axis x line=none,
	ymin=0, ymax=32,
	ytick={8,16,24,32}]
	grid=both,
	grid style={line width=.1pt, draw=gray!10},
	major grid style={line width=.2pt,draw=gray!500}
	\pgfplotsset{every axis legend/.append style={
			at={(1.05,-0.2)},
			anchor=north east}}
	\addplot [mark=triangle*,black,very thick,dash dot,mark size=3.5pt] table [x=t, y=LR, col sep=semicolon]  {Exp.csv};
    \addlegendentry{Rubber (length)};
    \addplot [mark=triangle*,red,very thick,dash dot,mark size=3.5pt] table [x=t, y=LF, col sep=semicolon]  {Exp.csv};
    \addlegendentry{Foamed rubber (length)};
	\end{axis}
	\end{tikzpicture}
	    \caption{Position of the diffusion front and length of the submerged rubber, respectively foam, as a function of time.}
    \label{Fig:Exp1}
\end{figure}
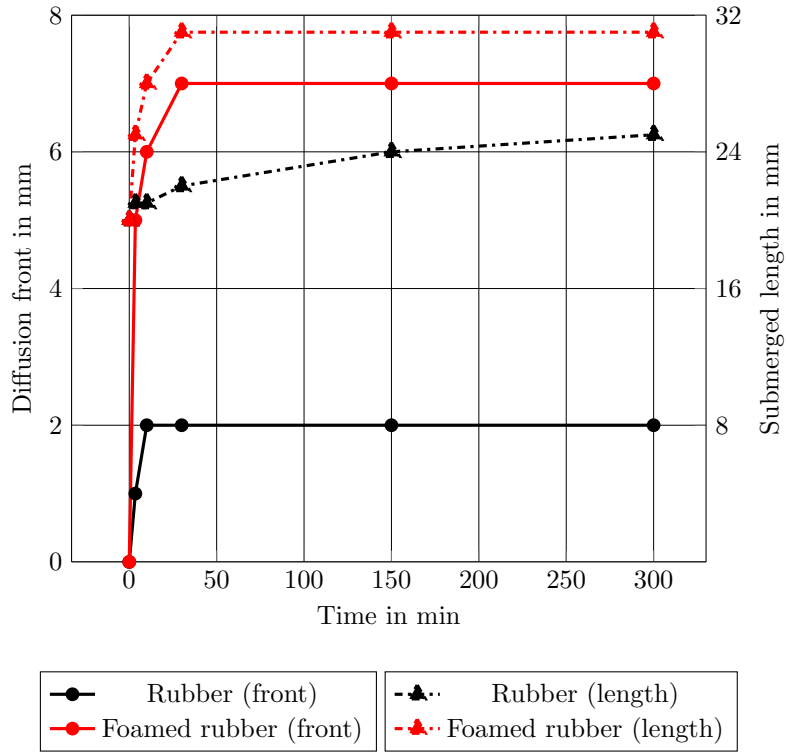

{\bf Experimental setup and results}:
A solution of cyclohexane and extracted black tea as marker is filled in a beaker. In order to compensate the evaporation of the solvent a dropping funnel is added, and adjusted accordingly to ensure a constant volume of solution inside the beaker, see Figure\,\ref{Fig:Exp}(left). The needed dropping rate is determined beforehand in dependency of the evaporation's rate. The test specimens are placed in the solution submerged by an initial length of 20 mm, see Figure\,\ref{Fig:Exp} (middle). To characterize the diffusion and swelling behaviour the diffusion front, the swelling length (submerged specimen) and the submerged swollen area are evaluated, see Figure\,\ref{Fig:Exp}(right), Figure\,\ref{Fig:Exp1}, Figure\,\ref{Fig:Exp2}, and  Table\,\ref{Tab:Exp}. For evaluation, the specimens are taken out of the solution, measured and placed back afterwards.\\

The position of the diffusion front outside the solution saturates for the rubber material before ten minutes and for the rubber foam before 30 minutes. The materials swell into the solution as well as it can be seen in Figure\,\ref{Fig:Exp2} which is indicated by the increasing submerged specimens lengths. In case of the rubber foam its saturation is reached in a similar time frame as for the diffusion front, whereas the full rubber only slowly converges. Due to the porosity of the foam and therefore lower (mechanical) resistance the effects are pronounced. Especially, the diffusion into the upper part of the specimen is amplified.  Similar but slower effects are observed by evaluating the swollen, submerged surface areas, see Figure\,\ref{Fig:Exp2}. 

\begin{table}[h]
	\centering
	\begin{tabular}{|c|c|c|c|c|c|c|}
		\hline
		 & \multicolumn{3}{c|}{Rubber} & \multicolumn{3}{c|}{Foamed rubber} \\
		\hline
		Time  & Front & Length & Area & Front & Length & Area \\ \hline
		[min] & [mm] & [mm] & [mm$^2$] &[mm] & [mm] & [mm$^2$] \\ \hline
		0 & 0 & 20 & 340 & 0 & 20 & 340 \\ \hline
		3.5 & 1 & 21 & 420 & 5 & 25 & 550 \\ \hline
		10 & 2 & 21 & 420 & 6 & 28 & 700 \\ \hline
		30 & 2 & 22 & 462 & 7 & 31 & 837 \\ \hline
		150 & 2 & 24 & 576 & 7 & 31 & 961 \\ \hline	
		300 & 2 & 25 & 625 & 7 & 31 & 961 \\ \hline
	\end{tabular}
	\caption{Determined experimental results related to diffusion front, submerged length and submerged area of the swollen specimens, cf.\,Figure\,\ref{Fig:Exp}.}
	\label{Tab:Exp}
\end{table}
\pgfplotsset{compat=1.8}
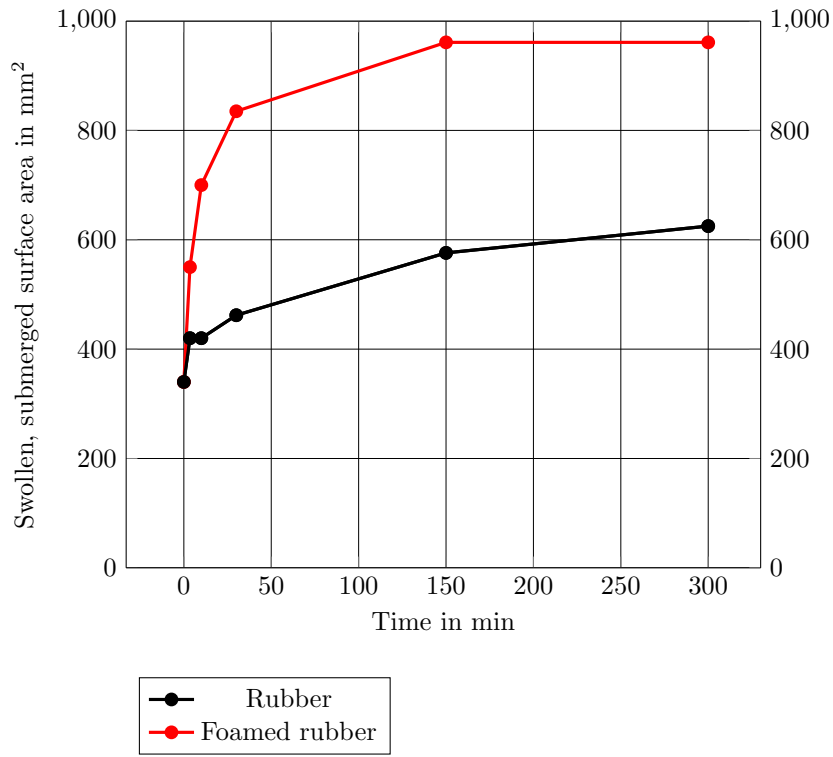
\begin{figure}[h!]
    \centering
   \begin{tikzpicture}
	\pgfplotsset{scale only axis,scaled x ticks=base 10:0,xmin=, xmax=330}
	\begin{axis}[
	axis y line*=left,
	ymin=0, ymax=1000,
	xlabel=Time in min,
	ylabel={Swollen, submerged surface area in mm$^2$},
	domain=0:350,
	grid=both,
	grid style={line width=.1pt, draw=gray!10},
	major grid style={line width=.2pt,draw=gray!500}]
	\pgfplotsset{every axis legend/.append style={
			at={(0.02,-0.2)},
			anchor=north west}}
	\addplot [mark=*,black,very thick] table [x=t, y=AR, col sep=semicolon]  {Exp.csv};
	\addlegendentry{Rubber};
	\addplot [mark=*,red,very thick] table [x=t, y=AF, col sep=semicolon]  {Exp.csv};
	\addlegendentry{Foamed rubber};
	\end{axis}
		\begin{axis}[
	axis x line=none,
	axis y line*=right,
	ymin=0, ymax=1000]
	\addplot [mark=*,black,very thick] table [x=t, y=AR, col sep=semicolon]  {Exp.csv};
\end{axis}
	\end{tikzpicture}
    \caption{Swollen, submerged surface area of the rubber, respectively foam, as a function of time.}
    \label{Fig:Exp2}
\end{figure}

\newpage

\section{Model equations}\label{model}

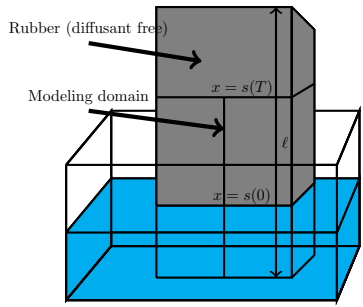
\begin{figure}[h]
	\begin{center}
		\begin{tikzpicture}[thick,scale=0.6, every node/.style={scale=0.60}]
		\draw (-2, -1.5) -- (4,-1.5) -- (4.5,-0.3) -- (-1,-0.3) -- (-2,-1.5); 
		\draw (-2, 0) -- (4,0) -- (4.5,1.2) -- (-1,1.2) -- (-2,0); 
		\draw (-2, 1.5) -- (4, 1.5) -- (4.5,2.7) -- (-1,2.7) -- (-2,1.5); 
		\draw (-2, 1.5) -- (-2, -1.5) -- (-1,-0.3) -- (-1,2.7) -- (-2,1.5); 
		\draw (4, 1.5) -- (4.5, 2.7) -- (4.5,-0.3) -- (4,-1.5) -- (4,1.5); 
		\filldraw[draw=black,fill=cyan] (-2, -1.5) -- (4,-1.5) -- (4.5,-0.3) -- (-1,-0.3) -- (-2,-1.5) ; 
		\filldraw[draw=black,fill=cyan] (-2, 0) -- (-2, -1.5) -- (-1,-0.3) -- (-1, 1.2) -- (-2,0); 
		\filldraw[draw=black,fill=cyan] (4, -1.5) -- (4.5, -0.3) -- (4.5,1.2) -- (4,0) -- (4,-1.5); 
		\filldraw[draw=black,fill=cyan] (-2, 0) -- (4,0) -- (4.5,1.2) -- (-1,1.2) -- (-2,0); 
		\filldraw[draw=black,fill=cyan] (-1, -0.3) -- (4,-0.3) -- (4,0) -- (-1,0) -- (-1,-0.3); 
		\draw (0,-1) -- (3,-1) -- (3,5) -- (0,5) -- (0,-1) ; 
		\draw (3,5) -- (3.5,4.5) -- (3.5,-0.5) -- (3,-1); 
		\draw (-2, 1.5) -- (-2, -1.5) -- (-1,-0.3) -- (-1,2.7) -- (-2,1.5); 
		\draw (4, 1.5) -- (4.5, 2.7) -- (4.5,-0.3) -- (4,-1.5) -- (4,1.5); 
		\draw (-2, -1.5) -- (4,-1.5) -- (4.5,-0.3) -- (-1,-0.3) -- (-2,-1.5); 
		\filldraw[draw=black,fill=cyan] (0, 0.6) -- (3,0.6) -- (3.5,1.1); 
		\filldraw[draw=black,fill=gray] (0, 0.6) -- (3,0.6) -- (3,5)--(0, 5)--(0,0.6);
		\filldraw[draw=black,fill=gray] (3,5)-- (3.5, 4.5)-- (3.5, 1.1)--(3,0.6) -- (3,5);
		\draw (0, 1.5) -- (4, 1.5);
		\draw (1.5,-1) -- (1.5,3); 
		\draw[<->] (2.65, -1) -- (2.65, 5) node[midway,anchor=west] {$\ell$} ;
		\draw (3, 5) -- (3, 0);
		\node at (1.9,0.8) {$ x= s(0)$};
		\draw (0, 3) -- (3,3) -- (3.5,3.3); 
		\node at (1.9,3.2) {$ x= s(T)$};
		\draw[ultra thick, <-]  (1,  3.8) -- (-1.5, 4.2)node[anchor=south] {Rubber (diffusant free)} ;
		\draw[ultra thick, <-]  (1.5,  2.3) -- (-1.5, 2.7)node[anchor=south] {Modeling domain} ;
		\end{tikzpicture}
		\caption{Sketch of the experimental setup. A rubber material is put in contact with a reservoir of diffusants.}
		\label{Fig:1}
	\end{center}
\end{figure}

At time $t = 0$, we consider a piece of a rubber denoted by $\Omega$ of vertical length $\ell>0$ placed in contact with a diffusant reservoir. The initial configuration of the rubber is divided by the interface $x= s(0)$ into two parts.  The lower part  $\Omega_1(0) = \Omega \cap \{x < s(0)\}$ is in contact with a solvent rich in diffusant concentration, while the upper part  $\Omega_2(0) =\Omega \cap \{x > s(0)\}$ is diffusant-free. Borrowing terminology from phase-change problems \cite{alexiades1992mathematical}, we start off from the assumption that the region  $\Omega_1(0)$, where the diffusant has already penetrated, is separated by  a sharp interface  from    the diffusant-free region $\Omega_2(0)$.  
To keep things simple, we look at a symmetric piece of material  so that we can consider the concentration of diffusants as a function of the one-dimensional space position $x$ and of the time variable $t\in(0,T)$.  Here $T$ is the observation time. 
As time elapses, the diffusants penetrate through the small pores of the rubber driven by absorption and diffusion processes. As a result, the initial position of the interface $x= s(0)$ changes continuously in time to new locations inwards the material. We denote the position of the moving diffusant penetration front by  $x=s(t)$.  It is worth observing that its location at time $t$ is {\em a priori} unknown and must be resolved together with determining the diffusant concentration inside the material.  The two regions of interest are  the diffusants-penetrated part and the diffusants-free part, viz.  $$\Omega_1(t) := \Omega \cap \{x < s(t): \ t\in(0,T)\} \;\;\text{(with diffusants)}$$   and  $$\Omega_2(t) :=  \Omega \cap \{x > s(t): \ t\in(0,T)\}\;\;\text{(without diffusants)},$$ respectively. We illustrate in Figure \ref{Fig:1} a sketch of the experimental setup we have in  mind. 

\begin{figure}[h]
	\begin{center}
		\begin{tikzpicture}[scale=1, every node/.style={scale=1}]
		\draw[|-|, color=red] (-1,-4) -- (6,-4);
		\draw[|-|, color=red] (6,-4) -- (9,-4);
		\draw (2, -4.1) -- (2, -3.9);
		\node at (-1,-4.4) {$ 0$};
		\node at (2,-4.4) {$s(0)$};
		\node at (6,-4.4) {$s(T)$};
		\node at (9,-4.4) {$\ell$};
		\draw[ultra thick, <-]  (1, -4) -- (1, -5)node[anchor=north] {Zone inside rubber with penetration of diffusants at $t=0$} ;
		\draw[ultra thick, <-]  (3.5, -4) -- (3.5, -3)node[anchor=south] {Diffusant-free rubber at $t=0$} ;
		\end{tikzpicture}
		\caption{Modeling domain.}
		\label{Fig:2}
	\end{center}
\end{figure}
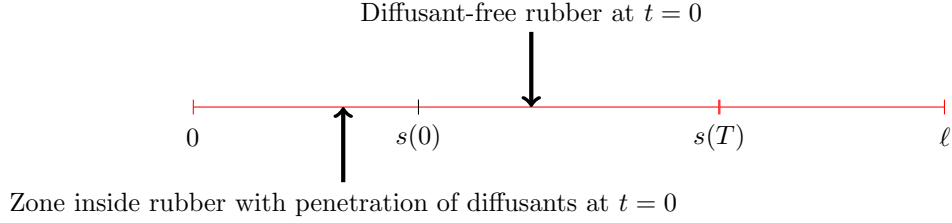

As $\Omega_2(t)$ is diffusants-free zone for each time $t$,  the actual  problem is to find the  diffusant concentration profile inside $\Omega_1(t)$ and the location of the free boundary $s(t)$.  This setting is called  one-phase moving boundary problem; see e.g. \cite{gupta2017classical} concerning modeling with moving interfaces. 
In this work, the modeling domain is a 1D slab as shown in Figure \ref{Fig:2}, which is the longitudinal line where  $0<s(0)\leq s(t)\leq \ell$.\\
For  a fixed observation time $T\in (0, \infty)$, the interval $[0, T]$ is the time span of the process we are considering.  Let  $x\in [0, s(t)]$ and $t \in [0, T]$ denote the space and respectively time variable,  and let $m(t, x)$  be the concentration of diffusant placed in position $x$ at time  $t$. 
The function $m(t, x)$ acts in the region $Q_s(T)$ defined by $$ Q_s(T):= \{ (t, x) | t \in (0, T) \; \text{and}\; x \in (0, s(t))\}.$$
The evolution of the diffusant concentration  is governed by the diffusion equation
\begin{equation} 
\label{a11} 
\displaystyle \frac{\partial m}{\partial t} -D \frac{\partial^2 m}{\partial x^2} = 0\;\;\; \ \text{in}\;\;\; Q_s(T),
\end{equation}
where $D>0$ is the constant diffusion coefficient. This choice indicates that we assume that the capillary suction does not  contribute to the transport of the diffusants in the rubber matrix; see \cite{Hewitt} for a case where the imbibition of a liquid into a solid matrix is driven by a suitable combination of viscous drag and capillary pressure. 
What concerns boundary conditions, we prescribe at the left boundary $x = 0$ an inflow. This is sometimes referred to as Newton or Robin-type boundary condition.  This is the place where absorption balances the incoming diffusion flux via 
\begin{equation}
\label{a12}
-D \frac{\partial m}{\partial x}(t, 0) = \beta(b(t) -Hm(t, 0))  \;\;\; \text{for}\;\; t\in(0, T). 
\end{equation}
In \eqref{a12}, $H>0$ is Henry's constant, while $\beta>0$ is a given absorption (or mass transfer) rate.  $\beta$ equals zero corresponds to the case where  there is no inflow of diffusants into rubber from the left boundary. Here we introduce the threshold $b$  to incorporate in the model an approximate value of the diffusant concentration that is initially present in contact with the rubber at the bottom surface as indicated in Figure \ref{Fig:1}.  The scenario in which diffusants tend to diffuse into  rubber corresponds to the case where the value of $b$ is greater than $Hm(0,0)$.  

At the right boundary, i.e. at $x = s(t)$, it holds: 
\begin{equation}
\label{a13}-D \frac{\partial m}{\partial x}(t, s(t))  =s^{\prime}(t)m(t, s(t))  \;\;\; \text{for}\;\; t\in(0, T).
\end{equation}
The boundary condition \eqref{a13} describes the mass conservation of diffusant concentration at the moving  boundary. It indicates that the diffusion mechanism is responsible for pushing the interface.\\
To have a complete model, we need to prescribe one additional condition on the right boundary as $s(t)$ is \textit{a priori} unknown. To this end, we propose an ordinary differential equation that describes explicitly the speed of the moving boundary by 
\begin{equation}
\label{a14}s^{\prime}(t) = a_0 (m(t, s(t)) - \sigma(s(t)) \;\;\;\;\text{for } \;\;\; t \in (0,T),
\end{equation}
where $a_0>0$ is constant. The  function $\sigma(s(t))$ incorporates an eventual contribution of the mechanics of the material. In particular \eqref{a14} points out that the mechanical behaviour (here it is about the swelling of the rubber) also contributes to the motion of the moving penetration front. It is worth mentioning that when the term $\sigma(s(t))$ vanishes, the model  equations  form the so-called one-phase  Stefan problem with  kinetic  condition. This approach was originally used to model the motion of the ice-water interfaces (compare e.g. \cite{alexiades1992mathematical}) as well as the motion of sharp carbonation reaction interfaces in concrete (see  \cite{aiki2011free,aiki2013large}).\\
The initial  diffusant concentration and the initial position of the moving front  are specified via
\begin{align}
\label{a15}&m(0, x) = m_0(x) \;\;\;\text{for}\;\;\; x \in [0, s(0)],\\
\label{a16} & s(0) = s_0>0\; \text{with}\;\; 0<  s_0< s(t) < L.
\end{align}
The model described in this Section was originally proposed in \cite{NHM} and analyzed mathematically in \cite{kumazaki2020global}. Now is the moment to verify its range of validity against the measurements reported in Section \ref{lab}.

\subsection{Non-dimensionlization}

To quantify the relative sizes of the characteristic time scales for absorption, diffusion and swelling, it is convenient to formulate the model equations in dimensionless form.  To do so, we introduce the dimensionless objects: 
\begin{align} 
\label{a17} z = \frac{x}{x_{ref}}, \hspace{7mm}  \tau = \frac{t}{x_{ref}^2/D}, \hspace{7mm}
u = \frac{m}{m_{ref}}, \hspace{7mm} h = \frac{s}{x_{ref}},
\end{align}
where  $x_{ref}$ is a characteristic  length scale,  while $m_{ref}$ is  a reference value for diffusant concentration.
The dimensionless concentration $u(\tau, z)$ acts in the region $Q_h(T^*)$  defined by $$ Q_h(T^*):= \left\{ (\tau, z) | \tau \in (0, T^*) \; \text{and}\; z \in (0, h(\tau)) \;\text{with}\; T^* = TD/x_{ref}^2 \right\}.$$
The   diffusion equation, now in dimensionless form, reads:
\begin{equation}
\label{a18} \frac{\partial u}{\partial \tau}  - \frac{\partial ^2 u}{\partial z^2} = 0 \;\;\; \ \text{in}\;\;\; Q_h(T^*).
\end{equation}
The boundary conditions \eqref{a12} and \eqref{a13} become
\begin{align}
\label{a19}& -\frac{\partial u}{\partial z}(\tau, 0) = \text{Bi} \left(\frac{b(\tau)}{m_{ref}} - Hu(\tau, 0)\right) \;\;\; \text{for}\;\; \tau\in(0, T^*), \\
\label{a20}&-\frac{\partial u}{\partial z}\left(\tau, h(\tau)\right) =h^{\prime}(\tau) u\left(\tau, h(\tau)\right)   \;\;\; \text{for}\;\; \tau\in(0, T^*),
\end{align}
where  Bi denotes the standard  mass transfer Biot number. In our context, this is defined as  \begin{equation} \label{Bi} \text{Bi} := \beta \frac{x_{ref}}{D}.\end{equation}
Equation \eqref{a14} takes the form
\begin{align}
\label{a21}h^{\prime}(\tau) = A_0 \left( u \left(\tau, h(\tau)\right) - \frac{\sigma(h(\tau))}{m_{ref}}\right)   \;\;\; \text{for}\;\; \tau\in(0, T^*),
\end{align}
where \begin{equation}\label{A0} A_0 :=  \frac{x_{ref}}{D}m_{ref} a_0. \end{equation} It is worth mentioning that $A_0$ is a sort of Thiele modulus (also called $2$nd Damk{\"o}hler number). Both dimensionless numbers Bi and $A_0$ are relevant for the discussion of our results as we will see in Section \ref{simulation}. The initial conditions become:
\begin{align}
\label{a22}&u(0, z) = \frac{m_0( z)}{ m_{ref}} \;\;\; \text{for}\;\;\; z\in[0,  h(0)],\\
\label{a23}& h(0) = \frac{s_0}{x_{ref}}.
\end{align}
Note also the following replacements 
\begin{align*}
&b(t_{ref} \tau) \rightarrow b(\tau),\\
&\sigma(x_{ref}h(t_{ref} \tau)) \rightarrow \sigma(h(\tau)),\\
&m_0(x_{ref} z) \rightarrow m_0(z).
\end{align*}
To fix ideas, we set as reference values \begin{align}
x_{ref} :=  s(T),\;\;\;\;\;
m_{ref} := m_0,
\end{align}
where $m_0$ is the initial amount of the solvent, while  the value $s(T)$ is the length of  penetration front when the solvent reached the stationary state. 

\subsection{Fixed-domain transformation}

To map the system  of equations \eqref{a18}-\eqref{a23} into the  cylindrical domain $Q(T^*) := \{ (\tau, y) | \tau \in (0, T^*) \; \text{and}\; y \in (0, 1)\}$, we use the transformation 
\begin{equation}
\label{a24}y = \frac{z }{h(\tau)},\;\;\; \tau \in (0, T^*).
\end{equation}
This transformation, originally due to Landau cf. \cite{landau1950heat}, maps the mathematical model with  moving boundary
conditions  into a mathematical model with fixed boundary conditions. By using the chain rule of differentiation and the transformation \eqref{a24}, equation \eqref{a18} becomes
\begin{equation}
\label{a32} 
\displaystyle\frac{\partial{u}}{\partial{\tau}} - y\frac{h^{\prime}(\tau)}{h(\tau)}  \frac{\partial{u}}{\partial{y}}- \frac{1}{(h(\tau))^2}\frac{\partial^2{u}}{\partial{y^2}} = 0
\;\;\; \ \text{in}\;\;\;Q(T^*),\\
\end{equation}
where we have used the prime ($^\prime$) notation to denote the derivative with
respect to time variable $\tau$.
The boundary conditions given in  \eqref{a19} and \eqref{a20} take now the form 
\begin{align} 
\label{a33} &- \frac{1}{h(\tau)} \frac{\partial u}{\partial y}(\tau, 0) = \text{Bi}\left(\frac{b(\tau)}{m_{0}}- Hu(\tau, 0)\right)\;\;\; \text{for}\;\; \tau\in(0, T^*), \\
\label{a34}&- \frac{1}{h(\tau)} \frac{\partial u}{\partial y}(\tau, 1) =  h^{\prime}(\tau)u(\tau, 1) \;\;\; \text{for}\;\; \tau\in(0, T^*), 
\end{align}
where $0<y<1.$
The ordinary differential equation corresponding to the  growth rate of the moving interface given in \eqref{a21} reads
\begin{equation}
\label{a35} h^{\prime}(\tau) =A_0\left(u(\tau, 1)-\frac{\sigma(h(\tau))}{m_{0}}\right)\;\;\; \text{for}\;\; \tau\in(0, T^*).
\end{equation}
The initial conditions  \eqref{a22} and \eqref{a23} transform into
\begin{align}
\label{a36} &u(0, y)  = 1 \;\;\; \text{for}\;\; y\in[0, 1], \\
\label{a37}&h(0) = \frac{s_0}{x_{ref}}.
\end{align}

\section{Numerical method}\label{FEM}

To handle  the discretization of $\eqref{a32}$ and \eqref{a35} we use the method of lines (see, for instance, \cite{thomee1984galerkin}). First, the equations are discretized in space by means of the finite element method. The resulting time-dependent ordinary differential equations (capturing information at the mesh points used for the finite element method)  are solved by the solver \texttt{odeint} in Python, see  \cite{johansson2018numerical} for details on the solver. In this Section, we present the details of our approximation scheme.\\
The weak formulation of the model is obtained by multiplying  \eqref{a32} with a test function $\varphi \in H^1(0, 1)$ and by integrating the result over the domain $(0, 1)$,
\begin{equation} 
\label{a38} 
\displaystyle \int_0^1\frac{\partial u}{\partial \tau} \varphi  dy- \int_0^1 y \frac{h^{\prime}(\tau)}{ h(\tau)}\frac{\partial u}{\partial y}\varphi  dy- \int_0^1\frac{1}{(h(\tau))^2}\frac{\partial^2{u}}{\partial{y^2}}\varphi dy  = 0.
\end{equation}
Using integration by parts  in \eqref{a38}  and  the boundary conditions  \eqref{a33} and \eqref{a34}, we get the following weak formulation: Find  the couple $(u, h)$ satisfying the equations 
\begin{align}
\nonumber&\displaystyle \int_0^1 \frac{\partial u}{\partial \tau} \varphi  dy- \int_0^1 y \frac{h^{\prime}(\tau)}{h(\tau)}\frac{\partial u}{\partial y} \varphi  dy +  \int_0^1\frac{1}{(h(\tau))^2}\frac{\partial{u}}{\partial{y}}\frac{\partial \varphi }{\partial y}  dy\\
\label{a39}&\hspace{4cm} - \frac{1}{h(\tau)} \text{Bi}\left(\frac{b(\tau)}{m_{0}}- Hu(\tau, 0)\right)\varphi(0) +\frac{h^{\prime}(\tau)}{h(\tau)}u(\tau, 1)\varphi (1) = 0 \;\;\text{for all} \; \varphi \in H^1(0,1),\\
& \label{a40} h^{\prime}(\tau) =A_0\left(u(\tau, 1)-\frac{\sigma(h(\tau))}{m_{0}}\right),\;\; \tau\in(0, T^*),
\end{align}
subject to the initial conditions \eqref{a36} and \eqref{a37}.

Existence and uniqueness of weak solution to problem~\eqref{a39}-\eqref{a40} are discussed in~\cite{kumazaki2020global}.  We refer the reader to Theorem 3.1 and Theorem 3.2 in~\cite{kumazaki2020global} for the details on these matters.  Furthermore, it turns out that the weak solution depends in a continuous way with respect to the model parameters. Now, we present approximation of this weak solution numerically.

We use here a Galerkin scheme based on piecewise linear polynomials to discretize the model equations \eqref{a39} and \eqref{a40}. For $N \in \{  2, 3,\cdots\}$, let $\{\phi_j\}_{j=0}^{N-1}$
denote the standard piecewise linear continuous basis function (i.e., the "hats  basis") for the finite dimensional subspace $V_k$ of $H^1(0, 1)$ defined on the uniform mesh $[0, 1/(N-1), 2/(N-1), \cdots, 1]$. The mesh size is $k := 1/(N-1)$. For $j\in \{ 0, 1,\cdots, N-1\}$, we have
\begin{equation*}
\phi_j(y) =  \begin{cases}
\displaystyle\;\;\frac{y - y_{j-1}}{k}, & \text{if\; $ y_{j-1} \leq y \leq y_j$},\\
\displaystyle\;\;\frac{y_{j+1}-y}{k}, & \text{if\; $  y_{j} \leq y \leq y_{j+1}$},\\
\;\;\;0, & \text{else}. \end{cases} 
\end{equation*}

\begin{figure}[h]
	\begin{center}
		\begin{tikzpicture}
		\draw [line width=0.35mm, ->](0,0) -- (0,4);
		\draw[line width=0.35mm,->] (0,0) -- (9,0);
		\draw[line width=0.25mm, dash dot] (0,2) -- (9,2);
		\draw (0,2) -- (2,0);
		\draw (0,0) -- (2,2) -- (4,0) ;
		\draw (2,0) -- (4,2) -- (6,0) ;
		\draw (4,0) -- (6,2) -- (8,0) ;
		\draw (6,0) -- (8,2)  ;
		\node at (0.2,2.2) {$ \phi_0$};
		\node at (2,2.2) {$ \phi_1$};
		\node at (4, 2.2) {$\phi_2$};
		\node at (6, 2.2) {$\phi_3$};
		\node at (8,2.2) {$\phi_4$};
		\node at (0,-.2) {$ 0=y_0$};
		\node at (2,-.2) {$ y_1$};
		\node at (4, -.2) {$y_2$};
		\node at (6, -.2) {$y_3$};
		\node at (8,-.2) {$y_4=1$};
		\node at (-0.2, 2) {$1$};
		\end{tikzpicture}
		\caption{Basis functions for $V_k$}
	\end{center}
\end{figure}
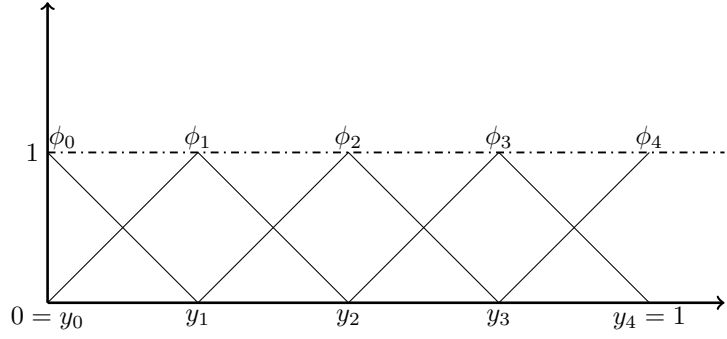

We formulate an approximate problem as follows:  Find the pair $(u_k, h)
\in C (0, T^*; V_k) \times C^1(0, T^*)$ such that
\begin{align}
\nonumber&\displaystyle \int_0^1 \frac{\partial u_k}{\partial \tau} \chi  dy- \int_0^1 y \frac{h^{\prime}(\tau)}{h(\tau)}\frac{\partial u_k}{\partial y} \chi dy +  \int_0^1\frac{1}{(h(\tau))^2}\frac{\partial{u_k}}{\partial{y}}\frac{\partial\chi }{\partial y}  dy\\
\label{a41}&\hspace{2cm}- \frac{1}{h(\tau)} \text{Bi}\left(\frac{b(\tau)}{m_{0}}- Hu_k(\tau, 0)\right) \chi (0) +\frac{h^{\prime}(\tau)}{h(\tau)}u_k(\tau, 1)\chi (1) = 0,\;\; \text{for all}\;\; \chi \in V_k,\;\;\tau\in(0, T^*) \\
\label{a42} & h^{\prime}(\tau) = A_0\left(u_k(\tau, 1)-\frac{\sigma(h(\tau))}{m_0}\right), \tau\in(0, T^*),\\
\label{a43}& u_k(0) = 1,\\
\label{a44}& h(0) = \frac{s_0}{x_{ref}}.
\end{align}
We let 
\begin{equation}
\label{a45}
u_k(\tau,y) = \sum_{j=0}^{N-1}\alpha_j(\tau) \phi_j(y),\;\; \tau\in [0, T^*] ,\;\; y\in[0, 1], 
\end{equation}
where
\begin{equation*} \alpha(\tau) = ( \alpha_0(\tau),  \alpha_1(\tau),  \cdots,  \alpha_{N-1}(\tau))^T\in \mathbb{R}^N. \end{equation*}
For $\tau = 0$,  $\alpha_j(0),\;\;j \in\{0,1,\cdots, N-1\}$ are the nodal values given by the initial approximation, i.e. for $j \in\{0,1,\cdots, N-1\}$, it holds:
\begin{equation} \label{a51}\alpha_j(0) =  1 .\end{equation}
Substituting \eqref{a45} in  \eqref{a41} yields
\begin{align}
\nonumber \sum_{j=0}^{N-1}\alpha_j^{\prime}(\tau)( \phi_j(y), \varphi) &-  \frac{h^{\prime}(\tau)}{h(\tau)} \sum_{j=0}^{N-1}\alpha_j(\tau)b( \phi_j(y), \varphi ) +\frac{1}{(h(\tau))^2} \sum_{j=0}^{N-1}\alpha_j(\tau)a( \phi_j(y), \varphi ) \\
\label{a46}& - \frac{1}{h(\tau)} \text{Bi}\left(\frac{b(\tau)}{m_{0}}- H\sum_{j=0}^{N-1}\alpha_j(\tau) \phi_j(0)\right) \varphi(0) \\
&+\frac{h^{\prime}(\tau)}{h(\tau)}\sum_{j=0}^{N-1}\alpha_j(\tau) \phi_j(1)\varphi(1) = 0,\;\; \tau\in(0, T^*). 
\end{align}
Taking, in particular, $\chi = \phi_i(y)\in V_k$, for $i \in\{0,1,\cdots, N-1\}$ in \eqref{a46},  we obtain

\begin{align}
\nonumber \sum_{j=0}^{N-1}\alpha_j^{\prime}(\tau)( \phi_j(y), \phi_i(y)) &-  \frac{h^{\prime}(\tau)}{h(\tau)} \sum_{j=0}^{N-1}\alpha_j(\tau)b( \phi_j(y), \phi_i(y) ) \\
\nonumber & +\frac{1}{(h(\tau))^2} \sum_{j=0}^{N-1}\alpha_j(\tau)a( \phi_j(y), \phi_i(y) ) 
- \frac{1}{h(\tau)} \text{Bi}\left(\frac{b(\tau)}{m_0}- H\sum_{j=0}^{N-1}\alpha_j(\tau) \phi_j(0)\right) \phi_i(0) \\
\label{a47}&   +\frac{h^{\prime}(\tau)}{h(\tau)}\sum_{j=0}^{N-1}\alpha_j(\tau) \phi_j(1)\phi_i(1) = 0,\;\; \tau\in(0, T^*) \;\;\;\text{for all}\;\;\;i \in\{0,1,\cdots, N-1\} .
\end{align}
We write this expression in matrix form as
\begin{align}
\nonumber M \alpha^{\prime}(\tau) &-\frac{h^{\prime}(\tau)}{h(\tau)}  K \alpha(\tau) + \frac{1}{(h(\tau))^2} A \alpha(\tau) \\
\label{a48}&-  \frac{1}{h(\tau)}  \text{Bi}\left(\frac{b(\tau)}{m_0}- H\alpha_0(\tau)\right) \textbf{e}_0 + \frac{h^{\prime}(\tau)}{h(\tau)}\alpha_{N-1}(\tau)\textbf{e}_{N-1}  = 0,\;\; \tau\in(0, T^*),
\end{align}
where  $\alpha(\tau) = (\alpha_0(\tau), \alpha_1(\tau), \cdots, \alpha_{N-1}(\tau))^T\in\mathbb{R}^N$ is the unknown vector. The matrices  $M, K, A \in \mathbb{R}^{N\times N} $  are defined by 
\begin{align*}
& M : = M_{i, j}= \int_0^1  \phi_i(y) \phi_j(y) dy,\;\;\; i,j \in \{0, 1, 2, \dots, N-1 \},\\
&M =
\frac{k}{6}\begin{bmatrix}
2 & 1 & 0 & 0 &0&  \cdots&0 &0  & 0&0 \\
1& 4 & 1 & 0 &0& \cdots &0 &0 & 0&0\\
0&1&4 & 1 & 0 & \cdots &0 &0 &0&0\\
\vdots&\vdots&\vdots &\vdots &  \vdots & \vdots&\vdots&\vdots&\vdots&\vdots\\
0&0&0&0&\cdots&0 & 1 & 4 & 1 &0 \\
0&0&0&0&\cdots&0 & 0 & 1 &4&1 \\
0&0&0&0&\cdots &0& 0 & 0 &1&2 \\
\end{bmatrix},\\
& K: = K_{i, j}= \int_0^1 y \phi_i(y) \phi_j^{\prime}(y) dy,\;\;\; i,j \in \{0, 1, 2, \dots, N-1 \},\\
& K = 
\begin{bmatrix}
\frac{2k- 3y_1}{6} & -\left( \frac{2k- 3y_1}{6} \right)& 0 &0&  \cdots&0 &0  & 0\\
-\left( \frac{2k- 3y_0}{6} \right)&\frac{4k - 3(y_0 + y_2)}{6} &-\left( \frac{2k- 3y_2}{6}\right)  &0& \cdots &0 &0 & 0\\
0&-\left( \frac{2k- 3y_1}{6} \right)&\frac{4k - 3(y_1 + y_3)}{6} &-\left( \frac{2k- 3y_3}{6}\right) & \cdots &0 &0 &0\\
\vdots&\vdots&\vdots &\vdots &  \vdots & \vdots&\vdots&\vdots\\
0&0&0&0&\cdots&-\left( \frac{2k- 3y_{N-3}}{6}\right)& \frac{4k - 3(y_{N-3} + y_{N-1})}{6} &-\left( \frac{2k- 3y_{N-1}}{6}\right) \\
0&0&0&0&\cdots &0 &-\left( \frac{2k- 3y_{N-2}}{6}\right)&\frac{2k - 3y_{N-2}}{6} \\
\end{bmatrix},\\
\vspace{2cm}
& A := A_{i, j}= \int_0^1 \nabla \phi_i(y) \nabla \phi_j(y) dy,\;\;\; i,j \in \{0, 1, 2, \dots, N-1 \},\\
&A = 
\frac{1}{k}\begin{bmatrix}
1 & -1 & 0 & 0 &0&  \cdots&0 &0  & 0&0 \\
-1&2 & -1 & 0 &0& \cdots &0 &0 & 0&0\\
0&-1&2 & -1 & 0 & \cdots &0 &0 &0&0\\
\vdots&\vdots&\vdots &\vdots &  \vdots & \vdots&\vdots&\vdots&\vdots&\vdots\\
0&0&0&0&\cdots&0 & -1 & 2 &-1&0 \\
0&0&0&0&\cdots&0 & 0 & -1 &2&-1 \\
0&0&0&0&\cdots &0& 0 & 0 &-1&1 \\
\end{bmatrix}.\\
\end{align*}
The vectors  $\textbf{e}_0, \textbf{e}_{N-1} \in \mathbb{R}^N$ are given by 
\begin{align*}
\textbf{e}_0 = (1, 0 , \cdots,0, 0)^T\; \text{and} \;\; \textbf{e}_{N-1} = (0, 0 , \cdots, 0, 1)^T.
\end{align*}
Note that the matrix $M$ is symmetric and positive definite as well as  diagonally dominant, and hence, invertible.
\begin{align}
\label{a49}
\nonumber&\alpha^{\prime}(\tau) = \frac{h^{\prime}(\tau)}{h(\tau)} M^{-1}K \alpha(\tau) -  \frac{1}{(h(\tau))^2} M^{-1} A \alpha(\tau)  \\
&\hspace{3cm} +   \frac{1}{h(\tau)}  \text{Bi}\left(\frac{b(\tau)}{m_0}- H\alpha_0(\tau) \right) M^{-1} \textbf{e}_1 - \frac{h^{\prime}(\tau)}{h(\tau)}\alpha_{N-1}(\tau)M^{-1}\textbf{e}_{N-1},\;\; \tau\in(0, T^*).
\end{align}
Substituting \eqref{a45} in \eqref{a42}, we get
\begin{align}
\label{a50}  &h^{\prime}(\tau) = A_0\left(\alpha_{N-1}(\tau)-\frac{\sigma(h(\tau))}{m_0}\right),\; \tau\in(0, T^*).
\end{align}
Equation \eqref{a49} together with \eqref{a50} forms a system of $N+1$ ordinary differential equations in terms of the unknowns
$(\alpha_0(\tau), \alpha_1(\tau), \cdots, \alpha_{N-1}(\tau), h(\tau))$. As the system is  not stiff, nor it contains algebraic constraints, it  can be solved numerically by  standard solvers such as  \texttt{odeint} in Python, relying on the initial condition $(\alpha_0(0),\alpha_1(0), \cdots, \alpha_{N-1}(0), h(0))$ given in  \eqref{a51} and \eqref{a44}. \\

\section{Model calibration. Simulation results for the dense rubber}\label{simulation}

In this section we present our simulation results for the dense rubber case. Basically, we calibrate here our model to represent the experimental range. The situation of the rubber foam is the aim of the next section.

The output consists of the concentration profile of the diffusants and of the position of the moving front. In this framework, we  investigate the parameter space by exploring eventual effects of the choice of parameters on the overall diffusants penetration process. Essentially, we deal with  \eqref{a49} and \eqref{a50} relying on repeated use of the solver \texttt{odeint} in Python. 
Our implementation requires three input arguments which are a vector of first order ordinary differential equations, a vector of initial conditions, and  a finite set of time discretization points. The output is a matrix where  each row is the solution vector of the ordinary differential equations defined at each discretization time point, producing the initial conditions vector in the first row of the output matrix. 
The simulation results are obtained using the set of reference parameters given in Table \ref{parameter}. 

\begin {table}[h]
\begin{center}
	\begin{tabular}{ |p{9.2cm}|p{1.7cm}|p{4.3cm}| }
		\hline
		Parameters & Dimension&Typical Values\\
		\hline
		Diffusion constant for concentration in rubber, $D$ &$L^2T^{-1}$& $3.66 \times 10^{-4}$ (mm$^2$/min), \cite{morton2013rubber} \\
		Absorption rate,  $\beta$ &$ LT^{-1}$ & $0.564$ (mm/min), \cite{rezk2018determination}\\
		Constant appearing in the speed of the moving boundary $a_0$&$L^4 T^{-1} M^{-1}$& $500$ (mm$^4$/min/gram) \\
		Initial height of diffusants (experiment),  $s_0$ & $L$ & $0.01$ (mm)\\
		$\sigma(s(t))$& $ML^{-3}$& $\frac{s(t)}{10}$ (gram/mm$^3$)\\
		Initial diffusant concentration (experiment), $m_0$ & $ML^{-3}$& $0.1$ (gram/mm$^3$)\\
		Concentration in lower surface of the rubber, $b$& $ML^{-3}$&$1$ (gram/mm$^3$) \\
		Henry's constant, 	$H$& --& 2.50 (dimensionless),  \cite{bohm1998moving}\\
		\hline
	\end{tabular}
	\caption {Name,  dimension and typical values for the model  parameters.}
	\label{parameter} 
\end{center}
\end {table}
Initially, the concentration is uniformly distributed within the rubber up to $0.01$ mm.  We take as observation time $T = 40$ minutes for the final time with time step $\Delta t = 1/1000$ minute. We choose the number of space discretization points $N$ to be $100$.  We take  the value $10$ mm  for characteristic length  scale  $x_{ref}$. As reference diffusant  concentration $m_{ref}$, we choose $0.1$ gram/mm$^3$. With our choice of parameters, the dimensionless numbers Bi and $A_0$ defined in \eqref{Bi} and \eqref{A0}, are of order of $O(10^4)$ and of  $O(10^6)$, respectively.

The  crucial component of the model is \eqref{a14} which describes the speed of the moving front. To get insight into  the role played by the parameters in  \eqref{a14}, we perform simulation for different values of $\sigma(s(t))$ and $a_0$. 

\subsection{Approximated swelling - capturing the effect of $\sigma(s(t))$}

\begin{figure}[h!]
	\centering
	\includegraphics[width=0.32\textwidth]{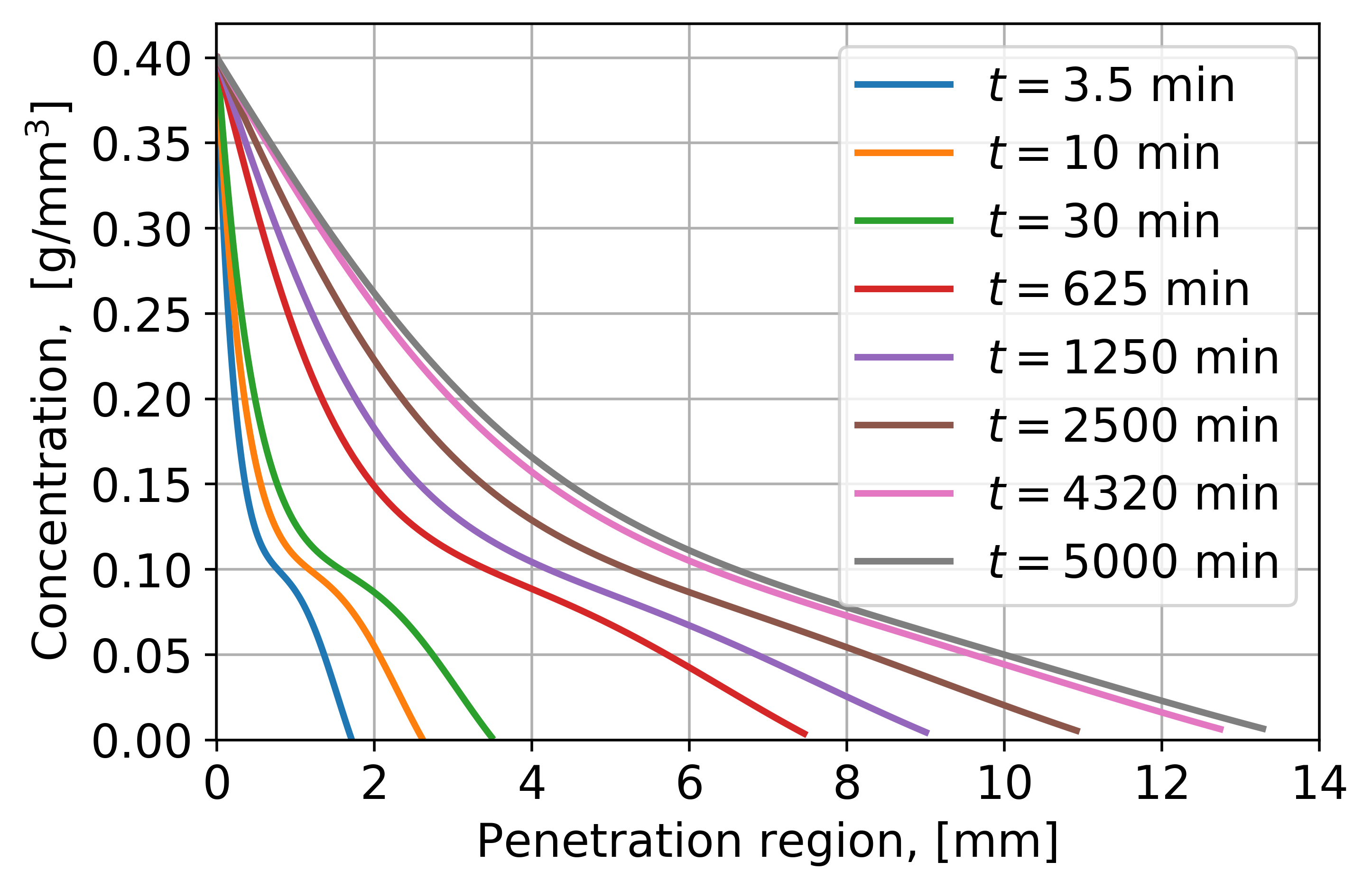}
	\hspace{0.1mm}
	\includegraphics[width=0.32\textwidth]{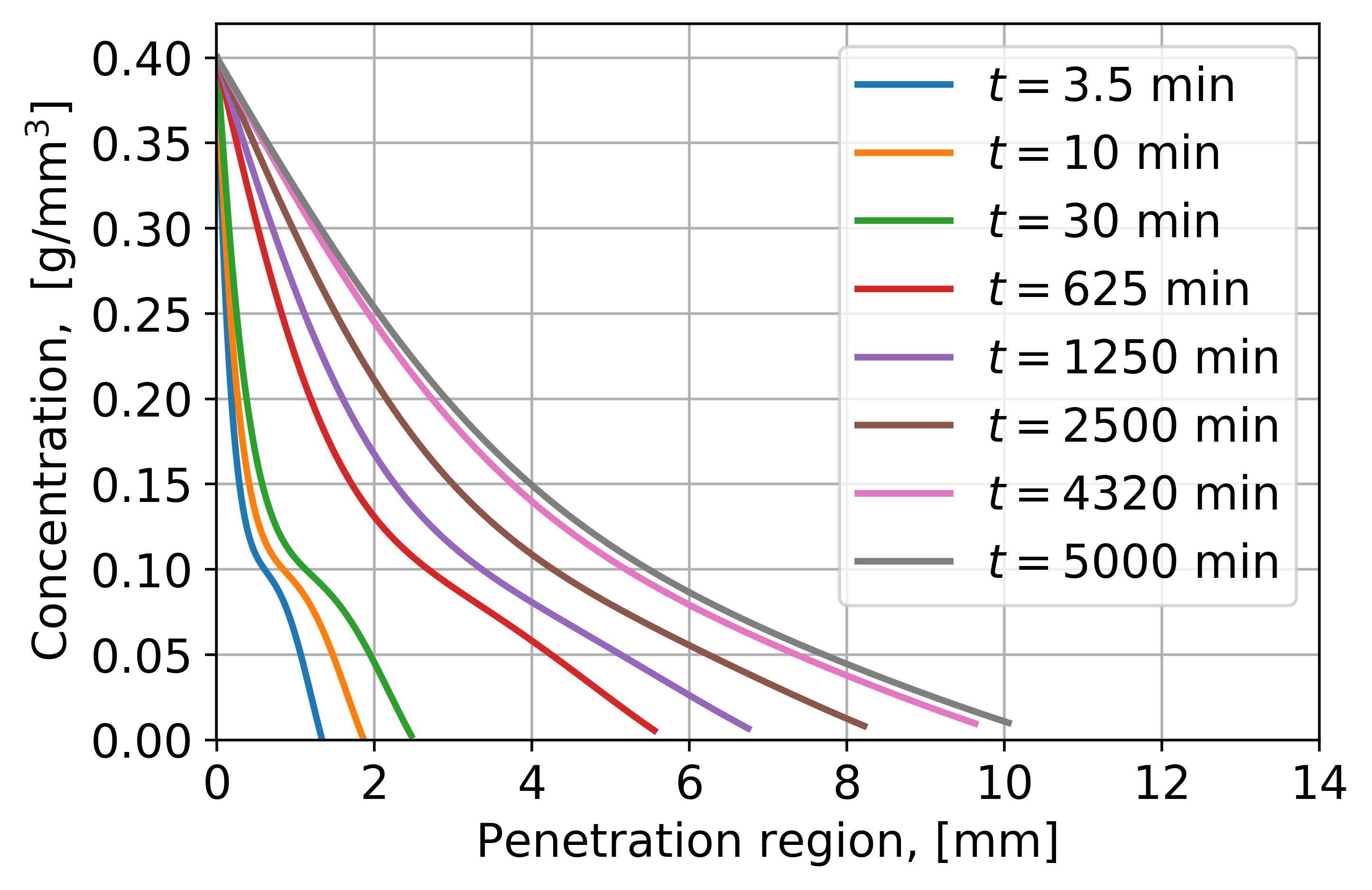}
	\hspace{0.1cm}
	\includegraphics[width=0.32\textwidth]{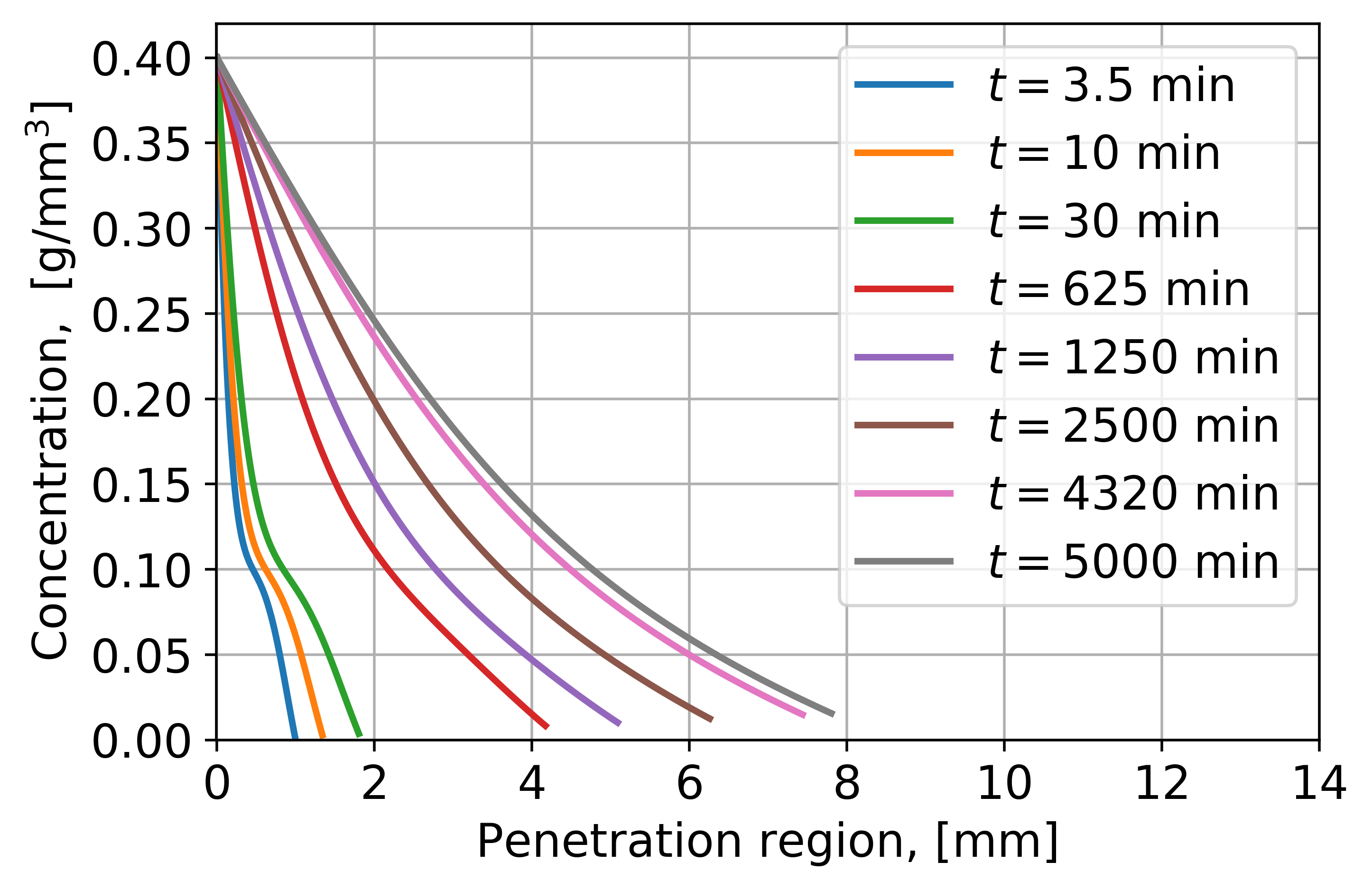}
	\caption{Concentration \textit{vs}. space with $\sigma(s(t)) =  \frac{s(t)}{20},\; \sigma(s(t)) = \frac{s(t)}{10}, \;\sigma(s(t)) = \frac{s(t)}{5}$ (from left to right).}
	\label{Fig:4}
\end{figure}

 We present numerical results for  different levels of linearity  of $\sigma(s(t))$.
 In Figure \ref{Fig:4} we  show the numerical results for the concentration profile for different values of $\sigma(s(t))$, while taking $a_0 = 500$. The first plot in Figure \ref{Fig:4}  shows a typical large-time behaviour very similar to the classical Stefan problem (ice-melting problem) where the concentration of the diffusant at the moving front position approaches zero.
Comparing now the remaining two plots in  Figure \ref{Fig:4}, we see that in this scenario the function $\sigma(s(t))$ describes a breaking mechanism preventing diffusants to move further inside the material. As a consequence of this fact, we notice an accumulation of diffusants at the position of the moving front.

It is worth noting that such an accumulation effect can be used in principle for material design purposes, since it gives an indication on how to bring a certain diffusant concentration in a prescribed region inside the material within a given time span. 

Considering  all three linear cases for $\sigma(s(t))$, within a short time
of release of diffusant from its initial position, the diffusant quickly enters the rubber from the left boundary (as we expect from \eqref{a13}) and then start off to diffuse further raising the penetration front. In the same Figure \ref{Fig:4}, we see that close to the initial time $t=0$, concentration profiles change convexity.  This situation disappears as time elapses. Such feature seems to be due to our choice of the parameter regime. In this particular case, a transient layer occurs in the time span of the process close to $t=0$. During this time span, there is apparently a rather strong  competition between diffusion, absorption, and swelling. The experimental data reported in Table \ref{Tab:Exp} refers to precisely this transient time during which most effects are visible. As expected, our numerical simulations indicate that the steady state of the process is not reached during the transient layer. 

\begin{figure}[h!]
	\centering
	\includegraphics[width=0.40\textwidth]{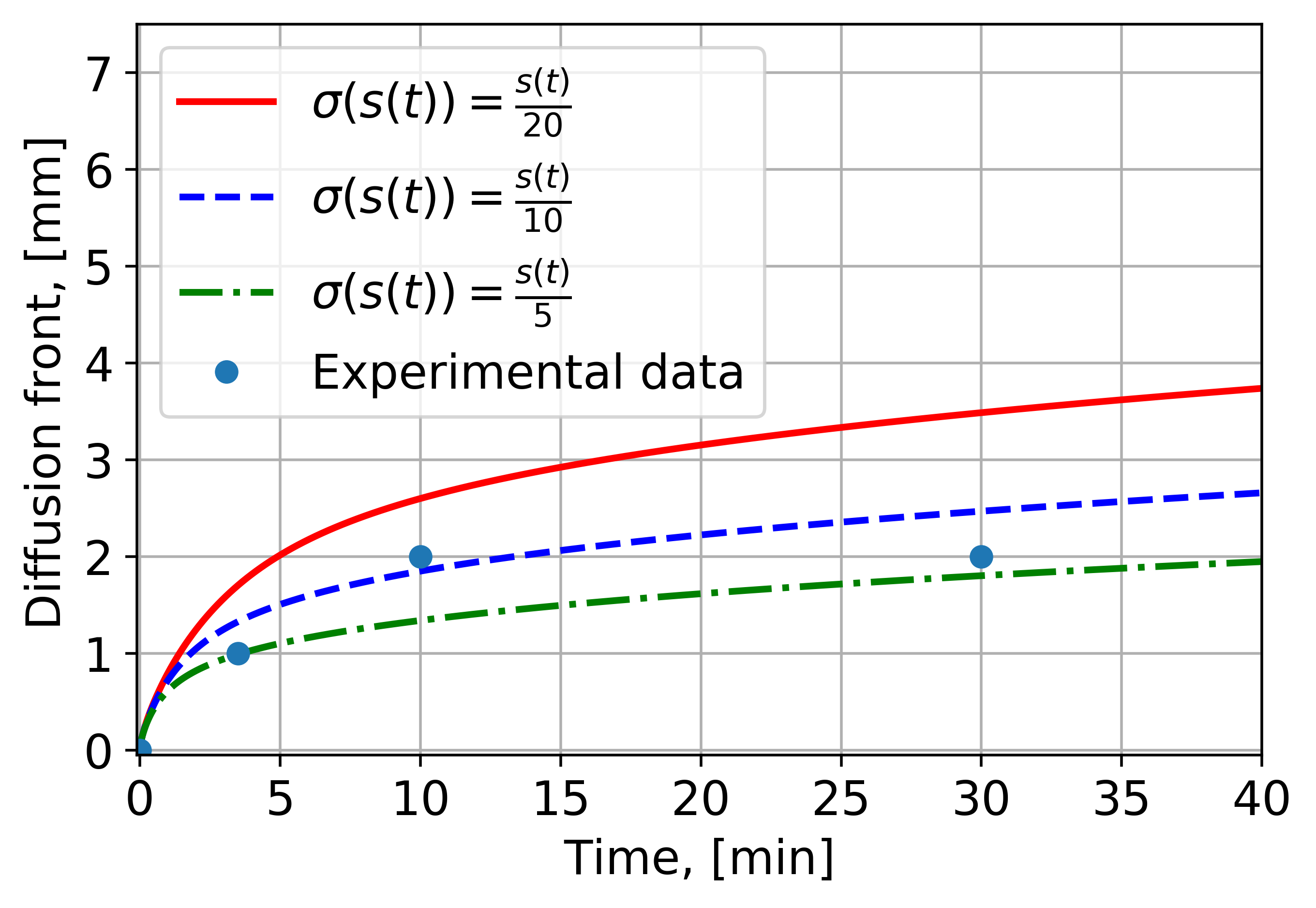}
	\caption{Comparison of experimental data and numerical diffusion front with different choices of  $\sigma(s(t))$ for $T= 40$ minutes.}
\label{Fig:10}
\end{figure}

In Figure \ref{Fig:10} we  compare the numerical results and experimental data for the moving front. The first and third plot in Figure \ref{Fig:10} show a large deviation between numerical results and  experimental data whereas the second plot shows a good agreement between numerical results and  experimental data. 

\subsection{Effect of the kinetic parameter $a_0$}

\begin{figure}[h!]
	\centering
	\includegraphics[width=0.40\textwidth]{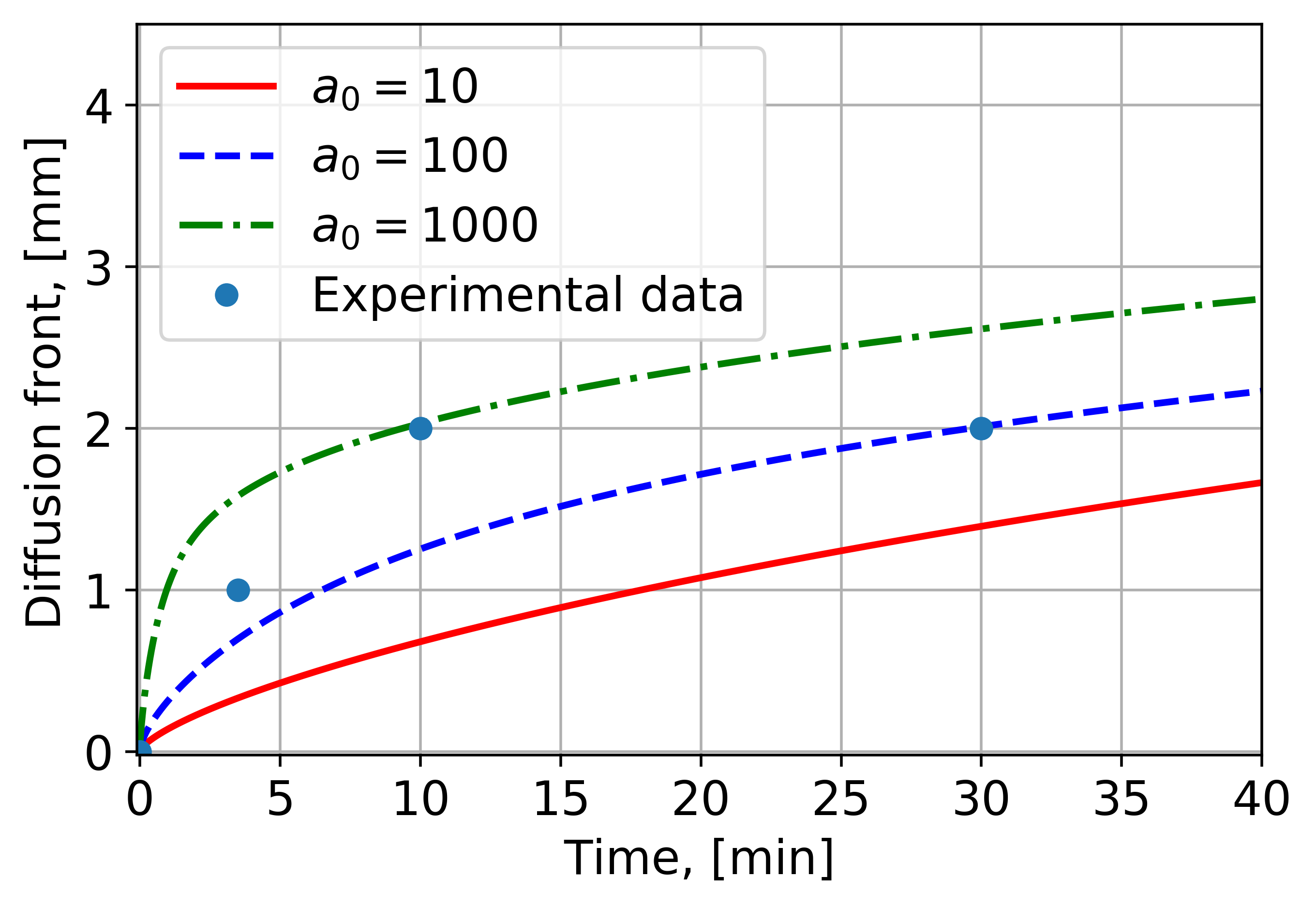}
\caption{Comparison of diffusion front for fixed $ \sigma(s(t)) = \frac{s(t)}{10}$ with varying $a_0 = 10,\; 100,\; 1000$.}
\label{Fig:6}	
\end{figure}

To get some insight in the role played by the parameter $a_0$,  we perform simulation runs for different values of $a_0$.  As illustrated in Figure \ref{Fig:6}, the changes in $a_0$ lead to  changes  in the numerical output. More precisely, larger values of $a_0$ lead to a faster penetration front, as expected in fact from \eqref{a14}.

\begin{figure}[h!]
	\centering
	\includegraphics[width=0.32\textwidth]{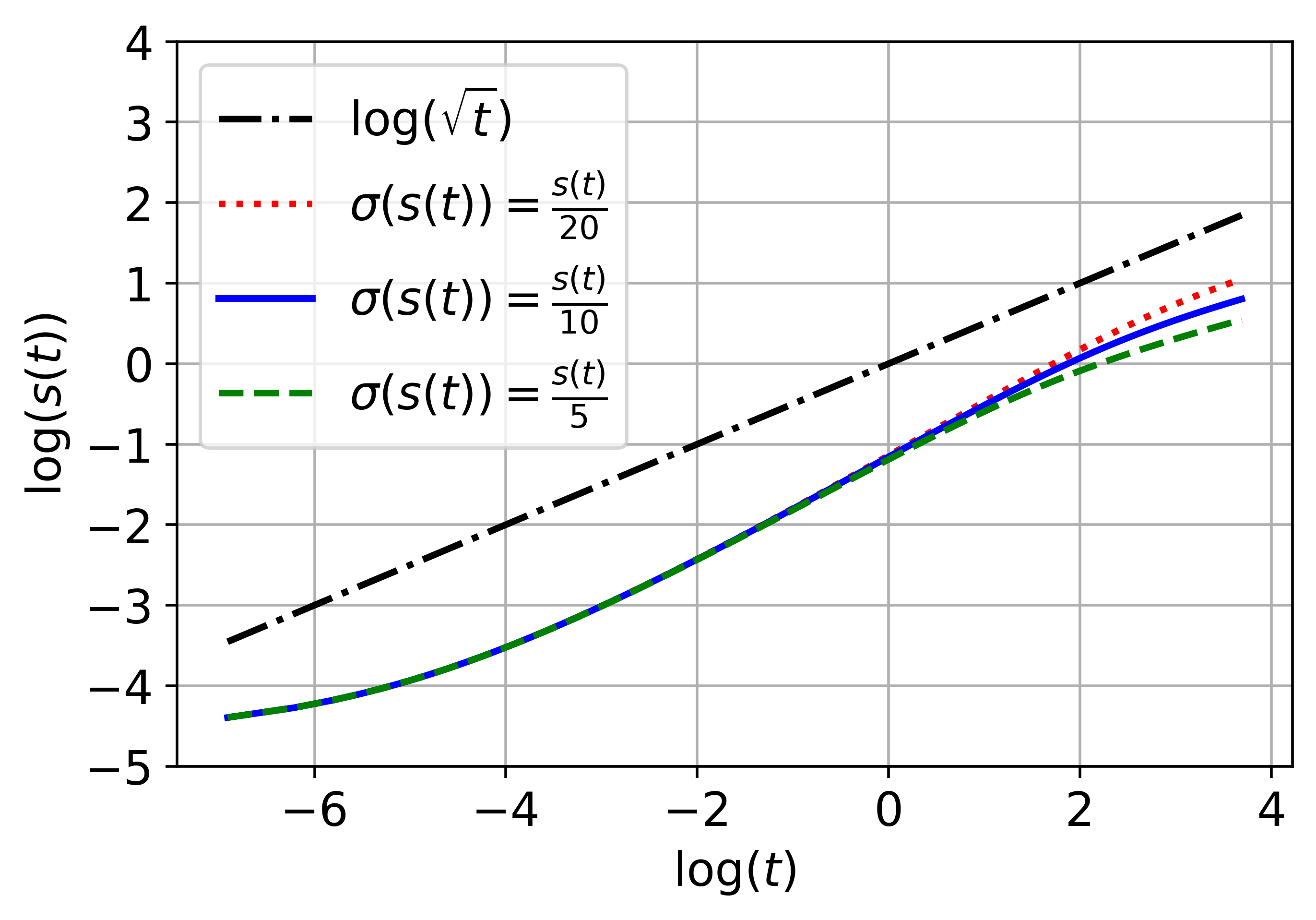}
	\hspace{0.1cm}
	\includegraphics[width=0.32\textwidth]{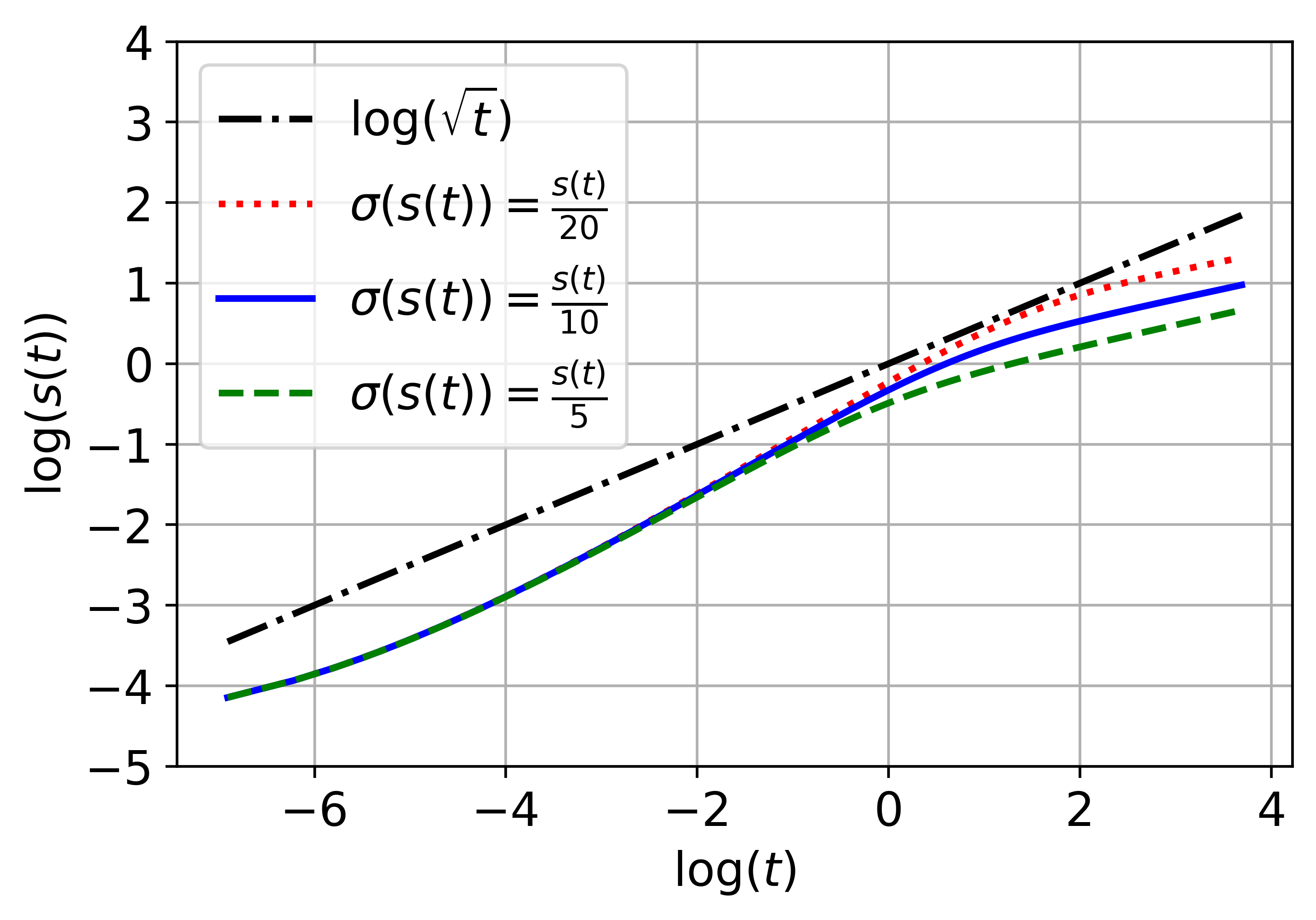}
	\hspace{0.1cm}
	\includegraphics[width=0.32\textwidth]{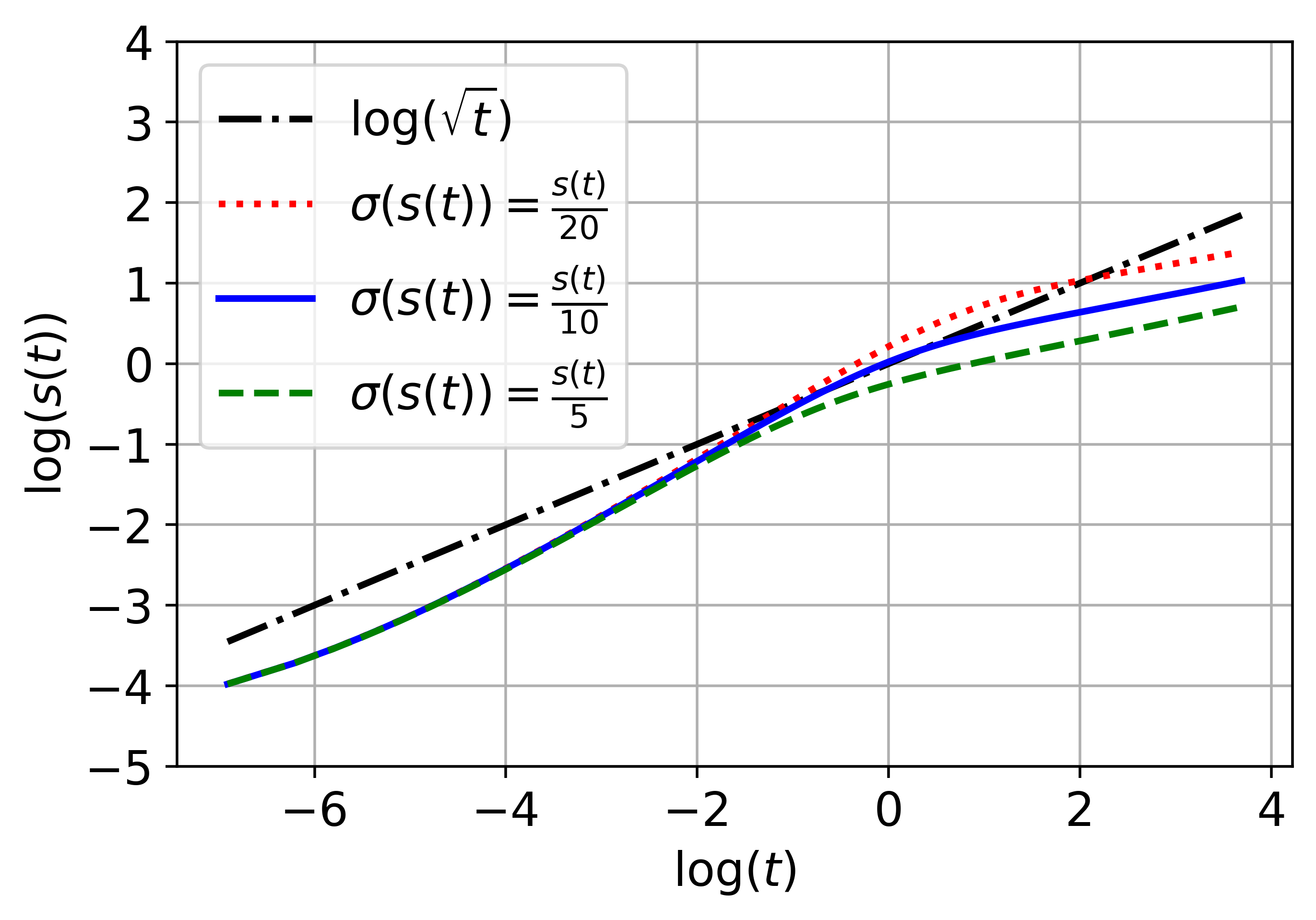}
	\caption{Short-time behaviour of  dense rubber: comparison of $\log(\sqrt{t})$ \textit{vs}. $\log(t)$ with  $\log(s(t))$ \textit{vs}. $\log(t)$ for different values of $\sigma(s(t))$ with $T = 40$ minutes and  $a_0 = 100, \;a_0 =  500, \;a_0 =  1000$ (from left to right).}
	\label{Fig:13}
\end{figure}

To estimate the dependence on time of the moving front $s(t)$, we fit our numerical results to the form 
\begin{equation} \label{a52}
s(t) =  t^{\gamma},\end{equation}
 where $\gamma> 0.$ In Figure \ref{Fig:13}, we show log-log plots of penetration fronts for different values of $\sigma(s(t))$ and $a_0$ and compare them against (\ref{a52}).    The slope of a log-log plot gives the value for $\gamma$, and a straight line is an indication that  relationship \eqref{a52} holds true. 
 The approximated values for $\gamma$ with varying the levels of linearity of $\sigma(s(t))$ and $a_0$ are listed in Table \ref{tab:1} and Table  \ref{tab:2}, respectively. We expect that the actual value of $\gamma$ strongly depends on the choice of structure of $\sigma(s(t))$.
 
 Considering now the information exhibited in Figure \ref{Fig:13}, we see comparisons of  the position of the front's penetration depth versus time for various modifications of the kinetic parameter $a_0$ as well as of the levels of linearity in the swelling term $\sigma(t)$. We compare all of them against the pure diffusive behavior $\sigma(t)\approx \sqrt{t}$. In other words, any deviation from a $\sqrt{t}$-like behavior can now be attributed to swelling.  Figure  \ref{Fig:13}  pinpoints what happens at the laboratory scale.
 Choosing the log-log scale in these plots allow us to identify the best exponents $\gamma$, which we list in Table \ref{tab:1}. At this stage, we can only say that if $\gamma> 0.5$ the motion of the diffusants penetration front appears to be super-diffusive, while for $\gamma< 0.5$ this is sub-diffusive.
	\begin {table}[h]
	\begin{center}
		\begin{tabular}{ |p{2.5cm}|p{2cm}|p{2cm}| p{2cm}|}
			\hline
			\diagbox{$a_0$}{$\sigma(s(t))$} & $\frac{s(t)}{20}$&$\frac{s(t)}{10}$& $\frac{s(t)}{5}$\\
			\hline
			100  &0.25684094& 0.18829908 & 0.11140718 \\
		500  &0.3726316 & 0.26559534 & 0.16501296\\
		1000&0.40093566& 0.28627561 &0.18102778 \\
			\hline
		\end{tabular}
		\caption {For the dense rubber: approximated value of  $\gamma$ entering $s(t) = t^\beta$ for $T = 40$ minutes.}
		\label{tab:1} 
	\end{center}
	\end {table}
	\newpage
\section{Simulation results for the foam rubber case}\label{Secfoam}

\begin{figure}[h!]
	\centering
	\includegraphics[width=0.32\textwidth]{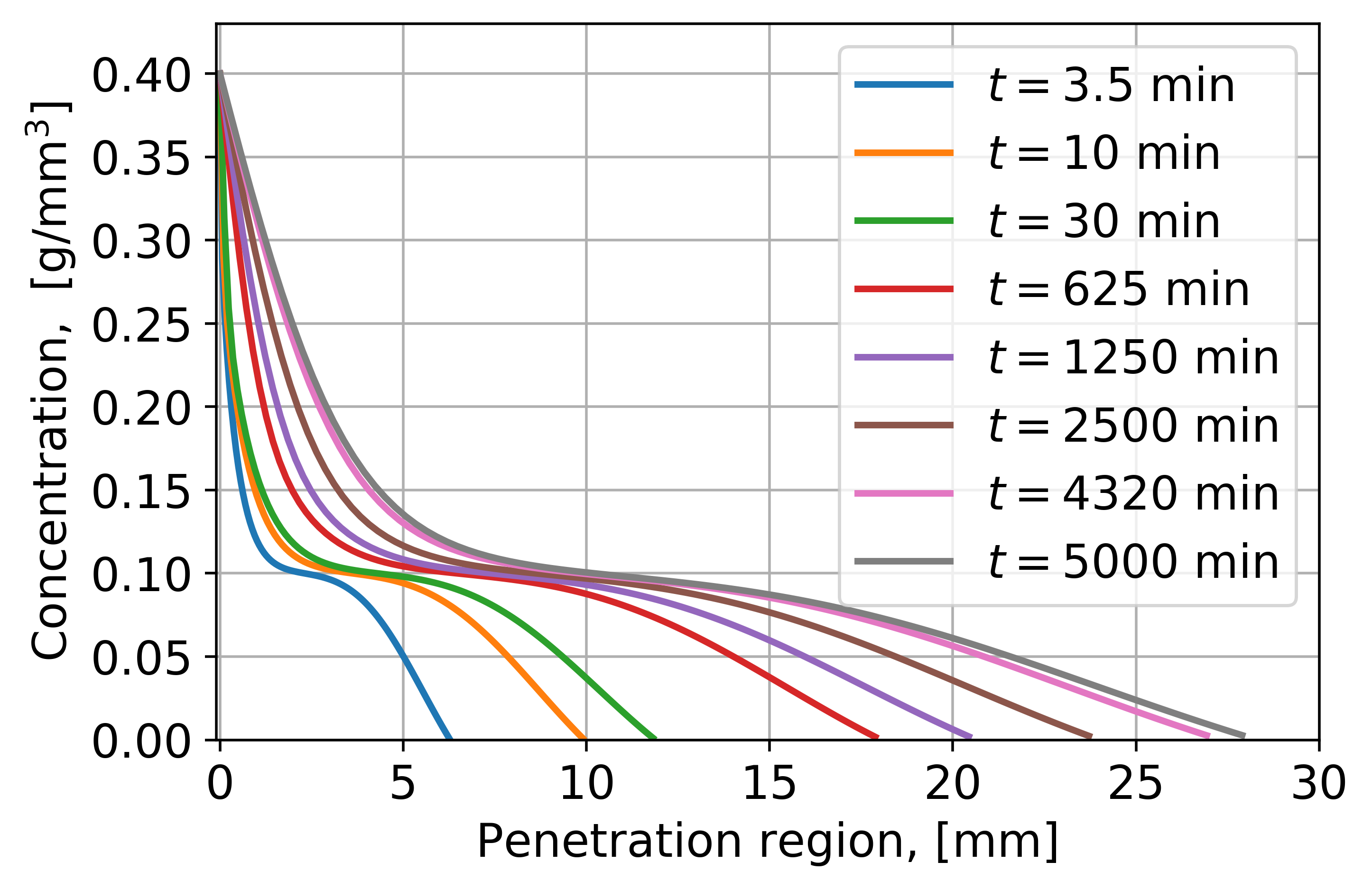}
	\hspace{0.1cm}
	\includegraphics[width=0.32\textwidth]{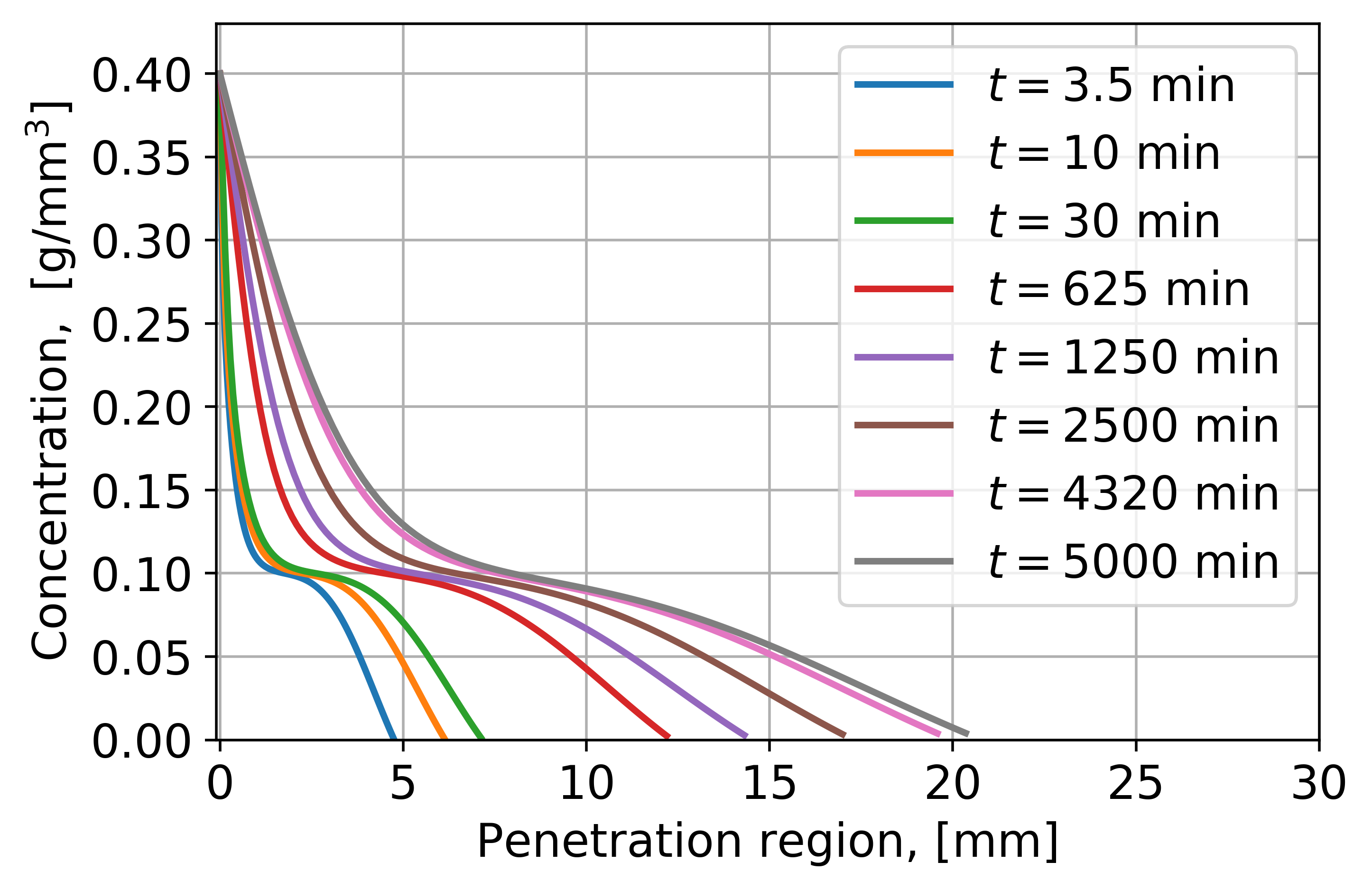}
	\hspace{0.1cm}
	\includegraphics[width=0.32\textwidth]{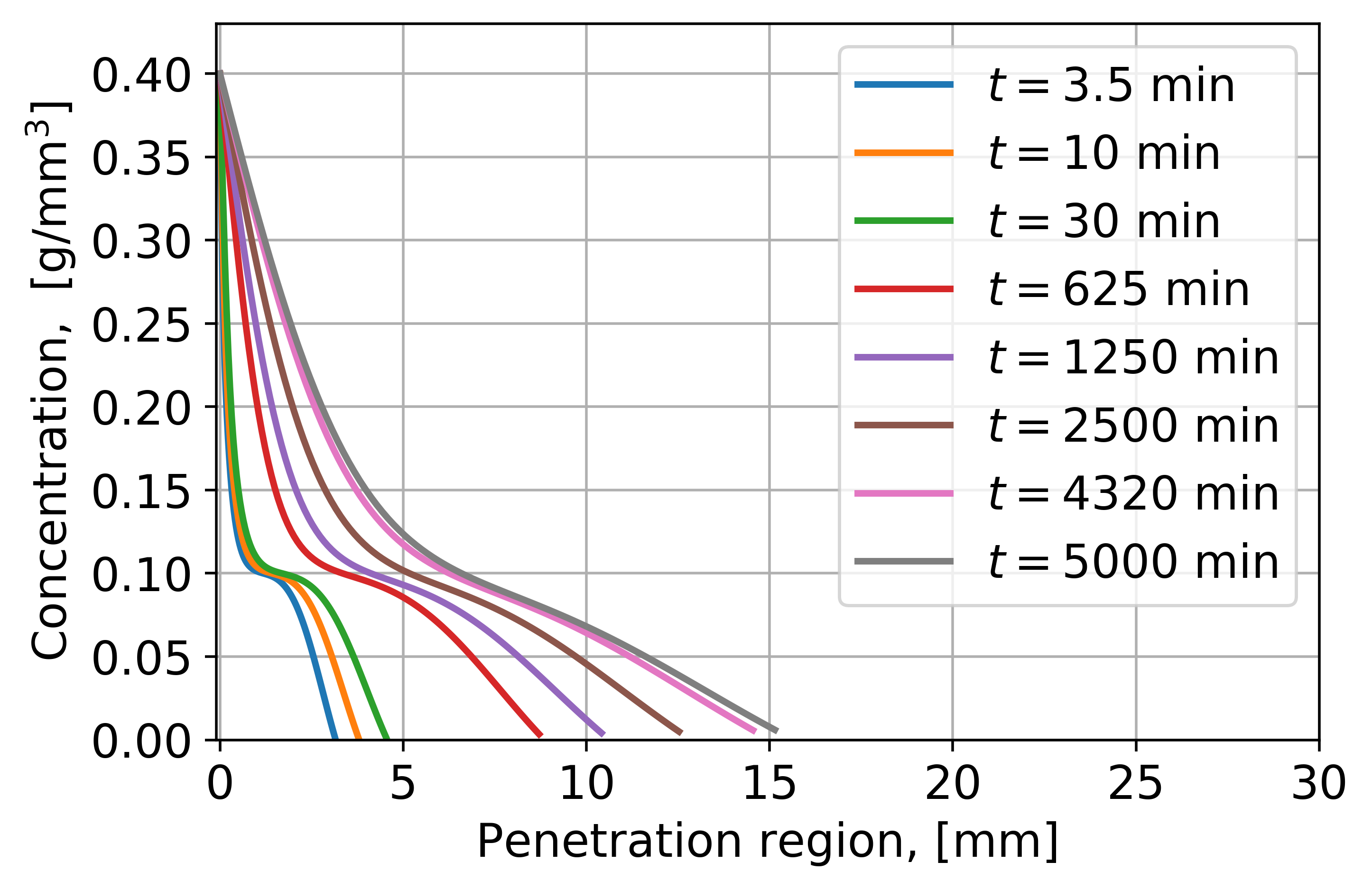}
	\caption{For the foam rubber case: concentration \textit{vs}. space with $\sigma(s(t)) =  \frac{s(t)}{100},\; \sigma(s(t)) = \frac{s(t)}{50}, \;\sigma(s(t)) = \frac{s(t)}{25}$ (from left to right).}
	\label{Fig:17}
\end{figure}

\begin{figure}[h!]
	\centering
	\includegraphics[width=0.40\textwidth]{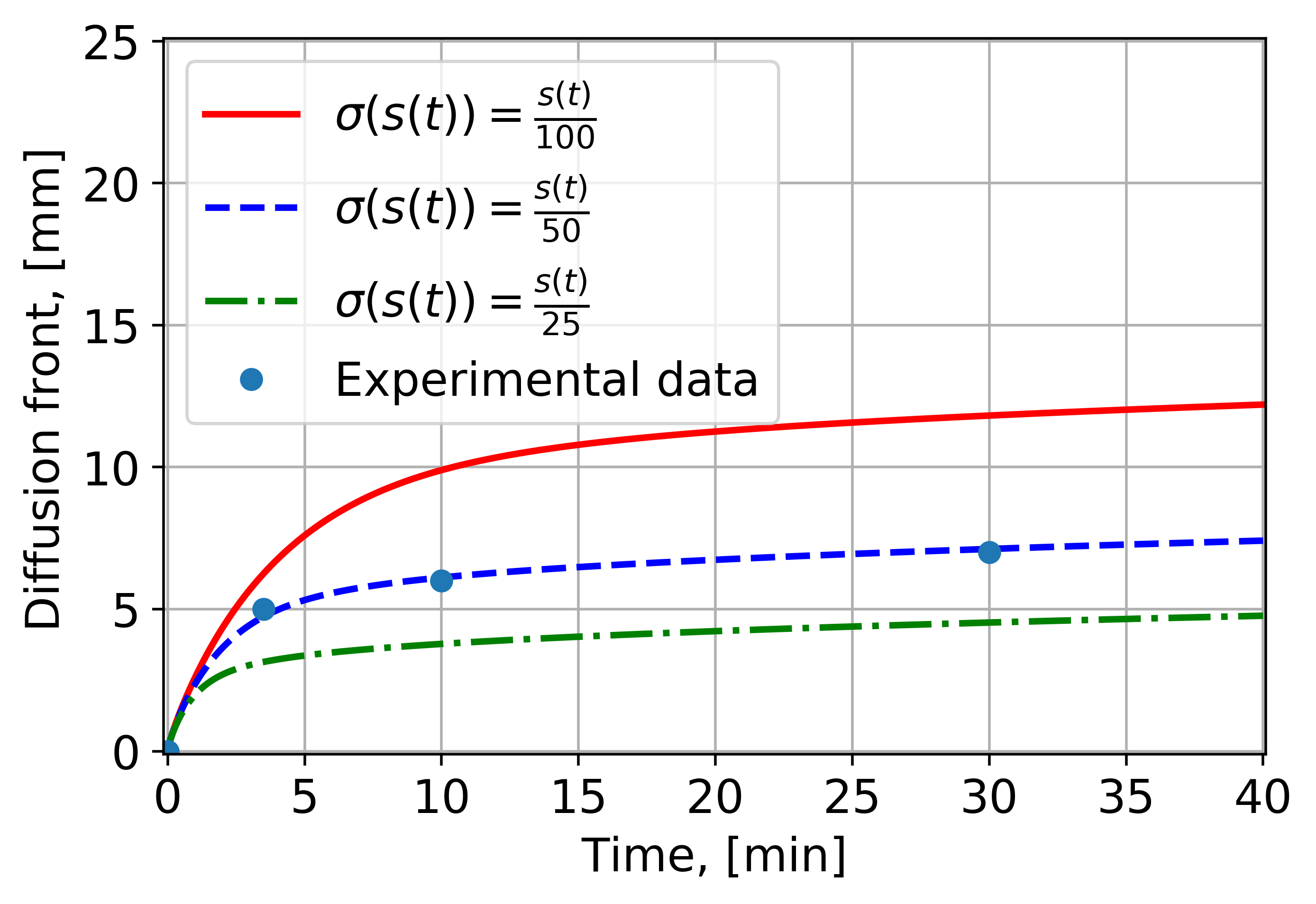}
	\caption{For the foam rubber case: comparison of the experimental diffusion front with numerical  diffusion front  for $a_0 = 2000$ and different choices of   $\sigma(s(t))$.}
	\label{Fig:20}
\end{figure}

\begin{figure}[h!]
	\centering
	\includegraphics[width=0.40\textwidth]{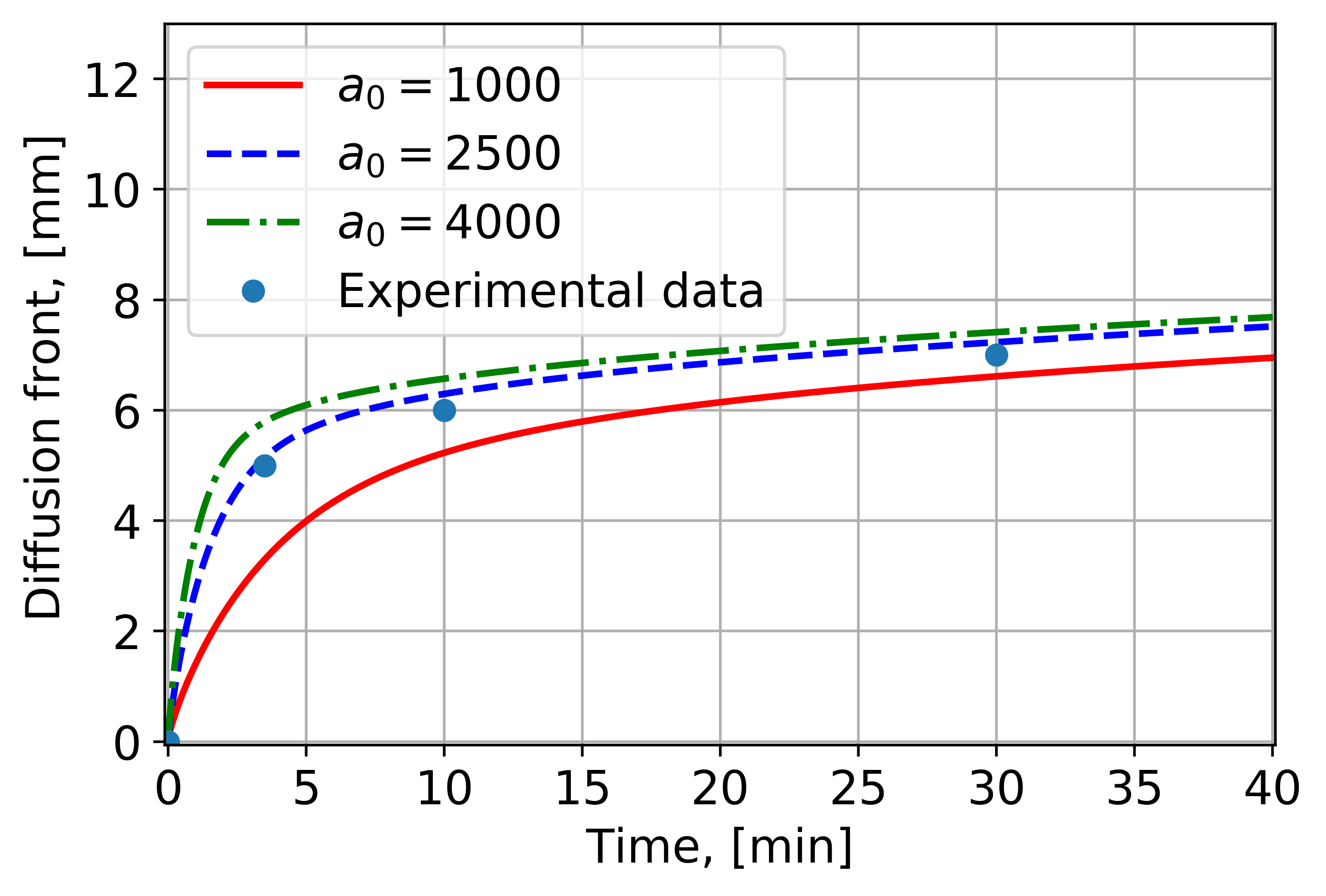}
	\caption{For the foam rubber case: comparison of the experimental diffusion front with  numerical  diffusion front  for $\sigma(s(t)) = \frac{s(t)}{50}$ and different choices of $a_0$.}
	\label{Fig:21}
\end{figure}

In this section, we perform numerical simulations to find out robust choices of  $\sigma(s(t))$ and $a_0$ to mimic the experimental data for the  foam rubber. Except for $\sigma(s(t))$ and $a_0$, we use  all other parameter values same as for dense rubber case discussed in Section \ref{simulation}. 
 In Figure \ref{Fig:17}, we show  the concentration profile for different values of $\sigma(s(t))$, while taking $a_0 = 2000$. With our choice of parameters, the dimensionless numbers Bi is in the same range as before, but $A_0$ is perhaps $1$ order of magnitude higher.
Comparing the plots in Figure \ref{Fig:17} with Figure \ref{Fig:4}, we see that the concentration profiles for the dense rubber and  foam rubber  cases have a similar shape.
In Figure \ref{Fig:20}, we compare  our numerical results with the experimental data for the moving front. One can observe some deviation between the numerical results and the experimental data in the case of the solid line (\full)  and dashdot (\dashdot), whereas the dashed line (\dashed) shows a good agreement. 
We also investigate the importance of the parameter $a_0$ what concerns  the case of  the  foam rubber. We compare the experimentally measured diffusion fronts against simulated diffusion fronts obtained when choosing different values of $a_0$ as shown in Figure \ref{Fig:21}.  Comparing the plots corresponding to $a_0 = 2500$ and $a_0 = 4000$, the effect of $a_0$ on the diffusion front is mainly seen during the first few minutes of the release of the concentration from its initial position, but not too much effect can be seen when time elapses beyond a couple of minutes.   

\begin{figure}[h!]
	\centering
	\includegraphics[width=0.32\textwidth]{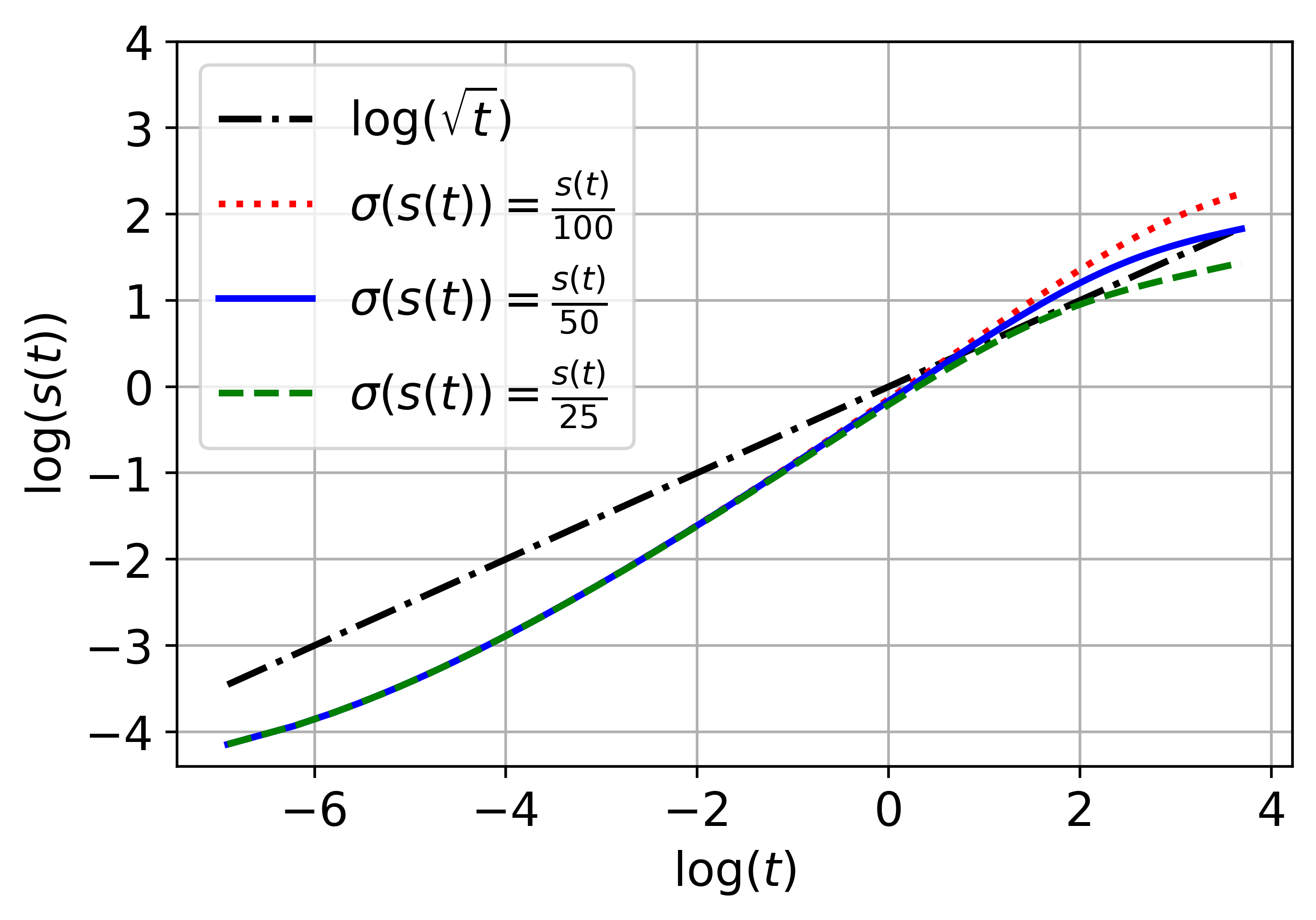}
	\hspace{0.1cm}
	\includegraphics[width=0.32\textwidth]{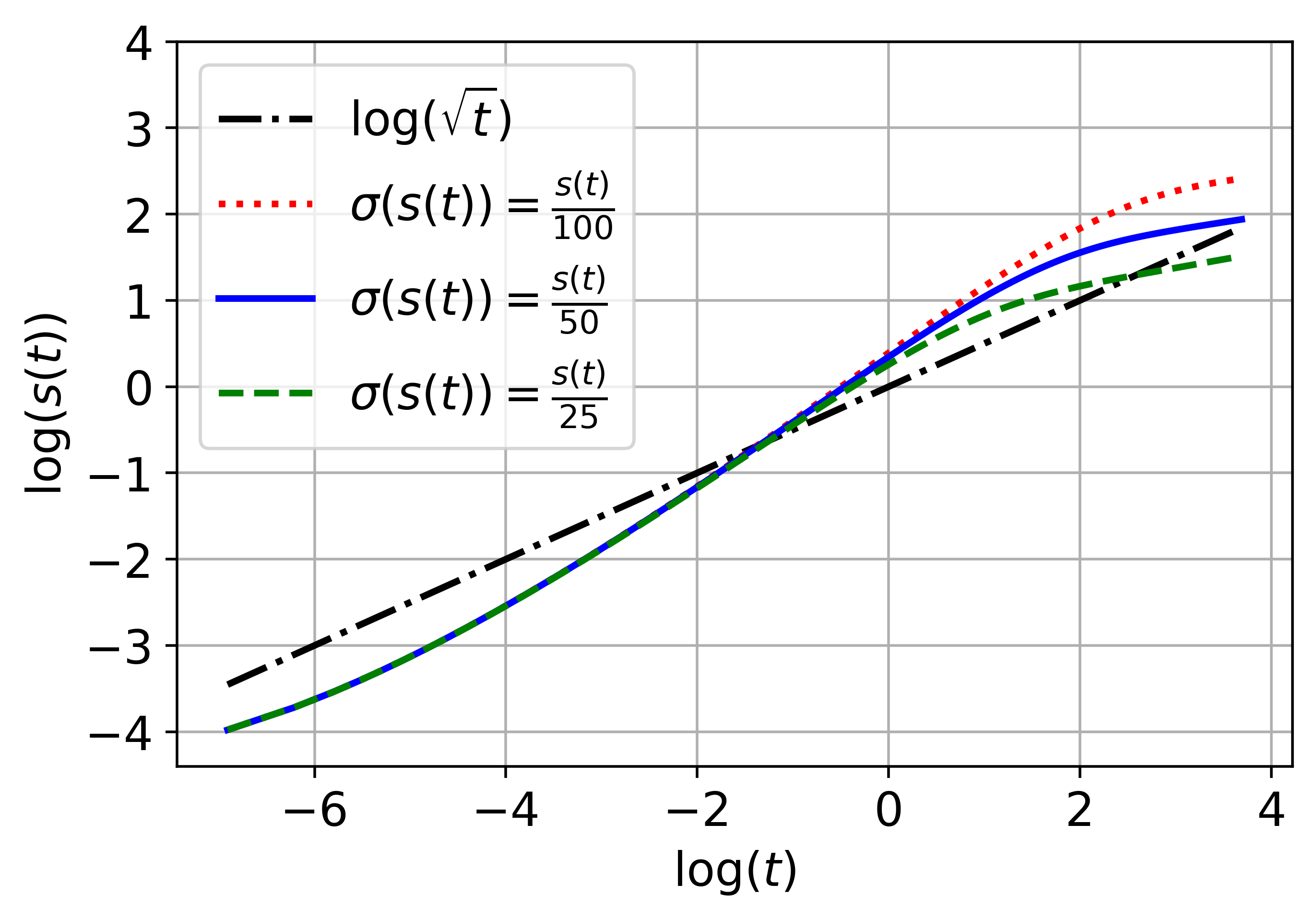}
	\hspace{0.1cm}
	\includegraphics[width=0.32\textwidth]{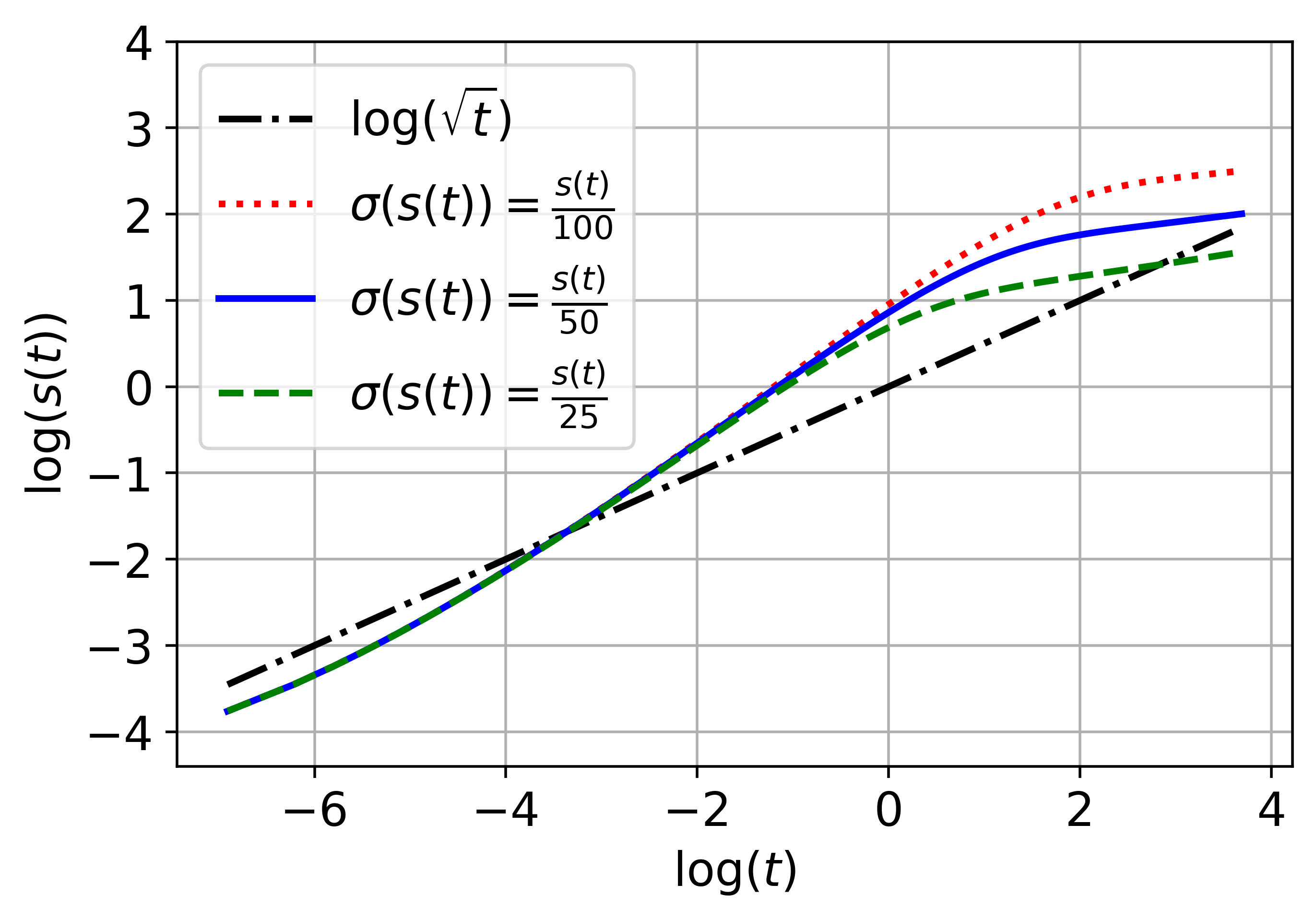}
	\caption{Short-time behaviour of  foam rubber : comparison of $\log(\sqrt{t})$ \textit{vs}. $\log(t)$ with  $\log(s(t))$ \textit{vs}. $\log(t)$ for different values of $\sigma(s(t))$ with $T = 40$ minutes and  $a_0 = 500, \;a_0 =  1000, \;a_0 =  2000$ (from left to right).}
	\label{Fig:22}
\end{figure}

We wish to see how much the asymptotic behavior of the diffusion fronts deviates from the $\sqrt{t}$-law.   To achieve this task, we  fit the simulated diffusion front  $s(t)$   to the form \eqref{a52}.
 We list in Table \ref{tab:3} and Table \ref{tab:4} the best fit values of $\gamma$ for different values of $\sigma(s(t))$ and $a_0$. Figure \ref{Fig:22} and Figure \ref{Fig:23} point out log-log plots of penetration fronts for different values of $\sigma(s(t))$ and $a_0$ compared against (\ref{a52}). 

	\begin {table}[h]
\begin{center}
	\begin{tabular}{ |p{2.5cm}|p{2cm}|p{2cm}| p{2cm}|}
		\hline
		\diagbox{$a_0$}{$\sigma(s(t))$} & $\frac{s(t)}{100}$&$\frac{s(t)}{50}$& $\frac{s(t)}{25}$\\
		\hline
		500 & 0.63059691 & 0.52306905 &  0.40861496\\
		1000  &0.70484987& 0.56983803 & 0.44029209 \\
		2000  &0.74367275 & 0.59585696  & 0.4592986 \\
		\hline
	\end{tabular}
	\caption {For the foam rubber case: approximated value of  $\gamma$ entering $s(t) = t^\gamma$ for $T = 40$ minutes.}
	\label{tab:3} 
\end{center}
\end {table}
\section{Expected large-time behavior of penetration fronts}\label{discussion}
We are interested in predicting the behavior of penetration depths of the diffusants beyond the laboratory timescale, hence the large-time asymptotic is now of interest. To perform such  large-time numerical simulations, we rely on the penetration depths listed in Table \ref{Tab:Exp} to which we add an additional point recorded after two days. This additional information was obtained in another large-time measurement performed for the same type of rubber with a similar solvent. The role of this last measurement point is to ensure that numerical results stay within the expected physical range as long as possible.

\begin{figure}[h!]
	\centering
	\includegraphics[width=0.32\textwidth]{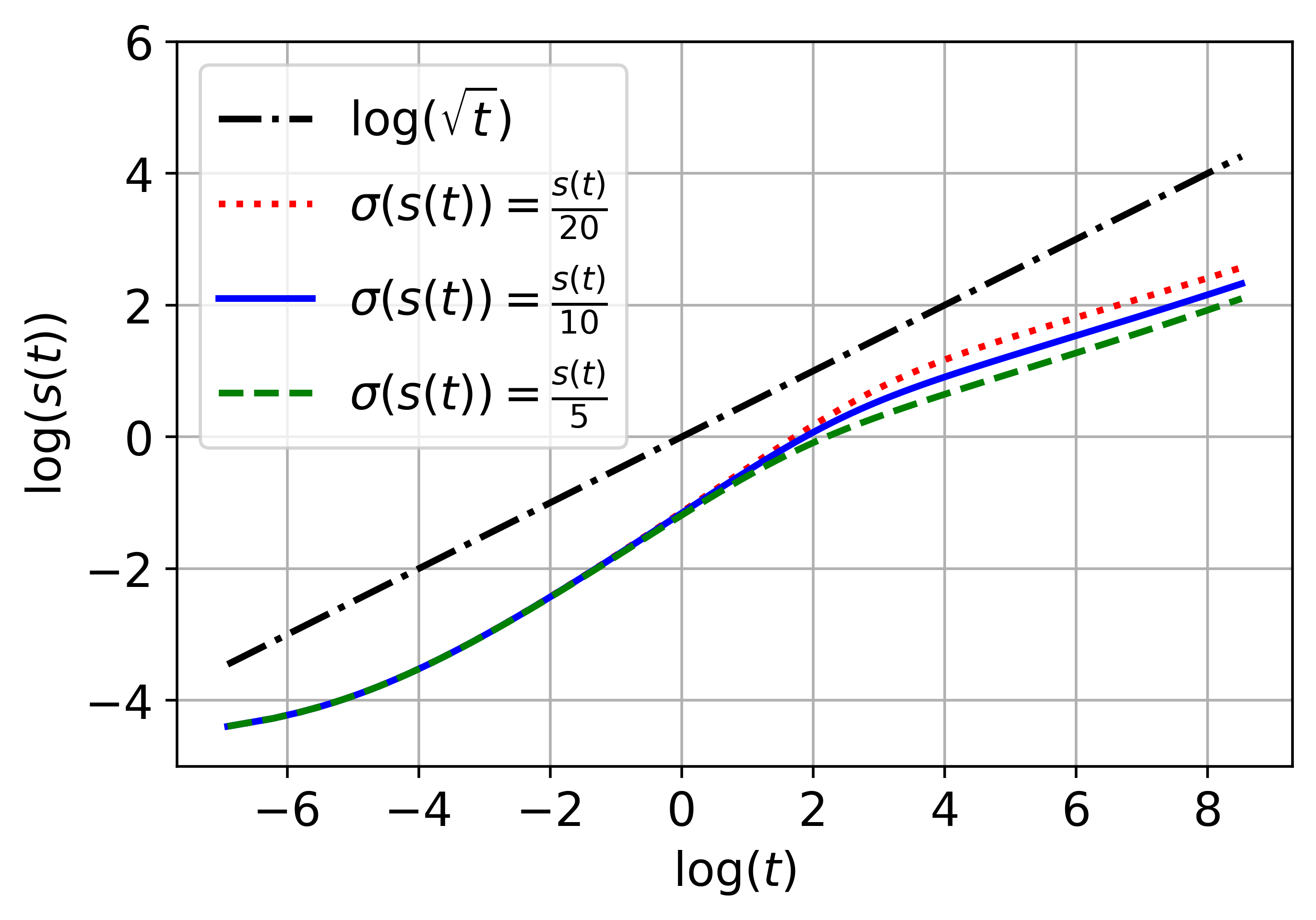}
	\hspace{0.1cm}
	\includegraphics[width=0.32\textwidth]{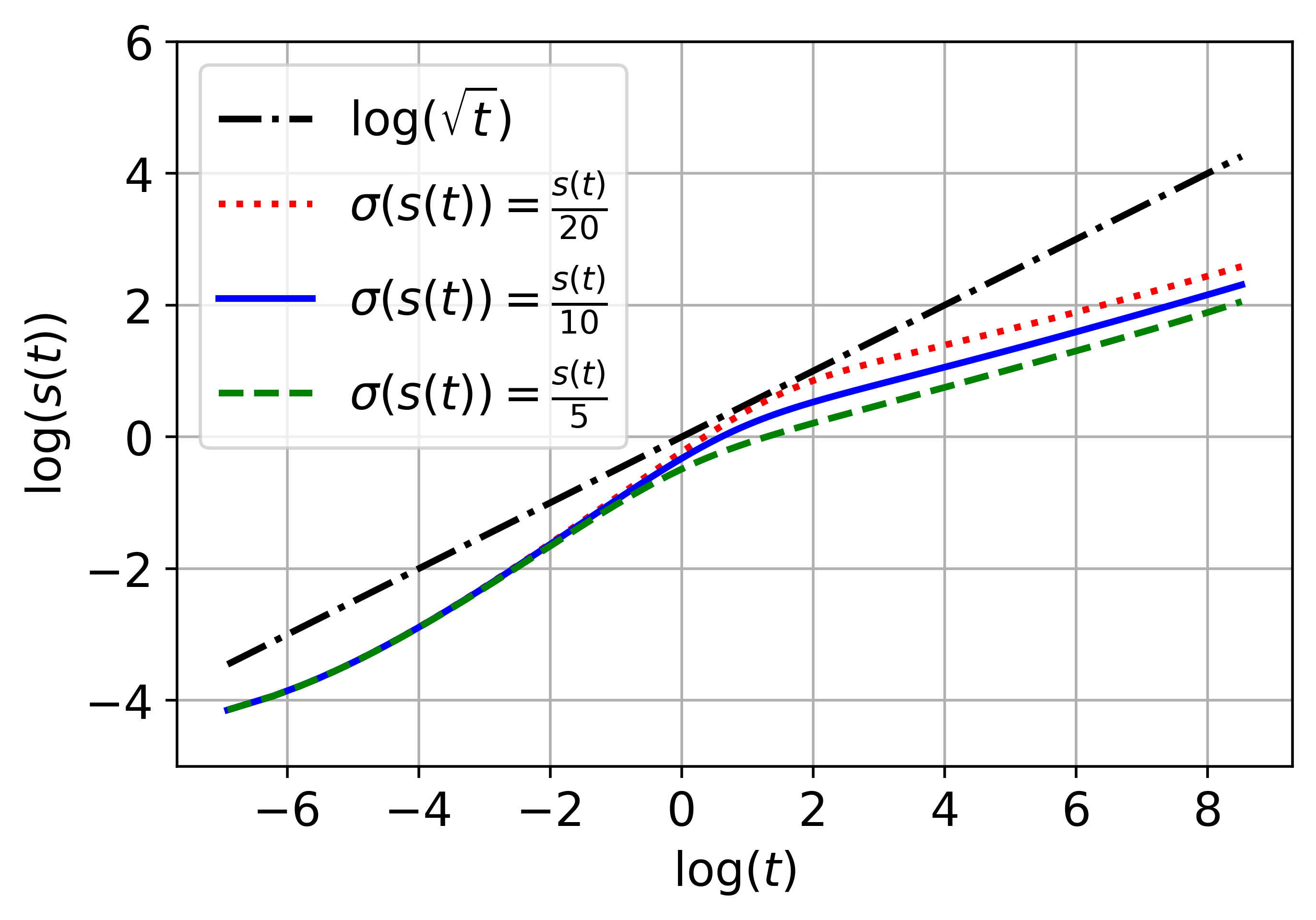}
	\hspace{0.1cm}
	\includegraphics[width=0.32\textwidth]{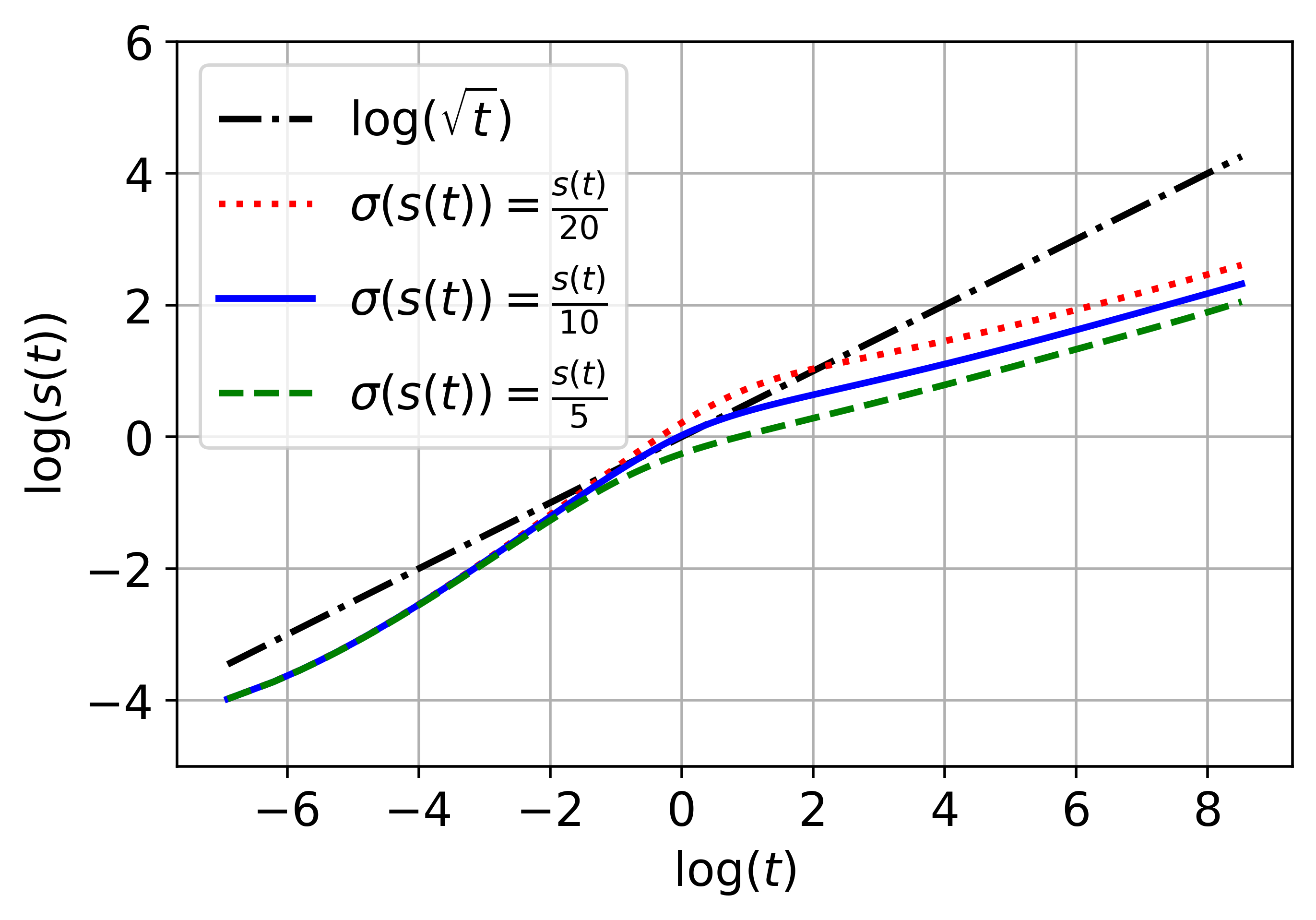}
	\caption{For the dense rubber case: comparison of $\log(\sqrt{t})$ \textit{vs}. $\log(t)$ with  $\log(s(t))$ \textit{vs}. $\log(t)$ for different values of $\sigma(s(t))$ with $T = 5000$ minutes and $a_0 = 100, \;a_0 =  500, \;a_0 =  1000$ (from left to right).}
	\label{Fig:15}
\end{figure}

\begin{figure}[h!]
	\centering
	\includegraphics[width=0.32\textwidth]{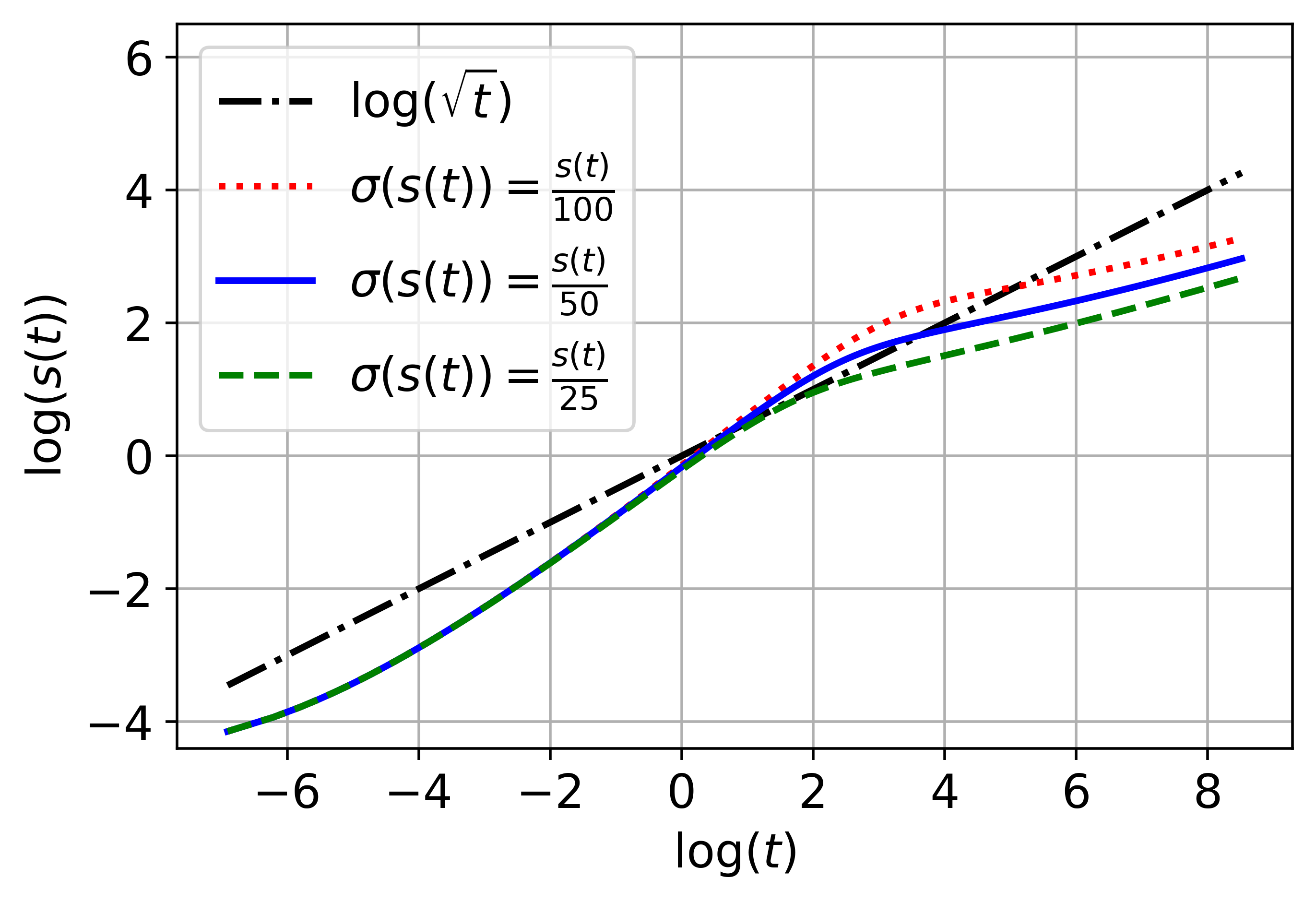}
	\hspace{0.1cm}
	\includegraphics[width=0.32\textwidth]{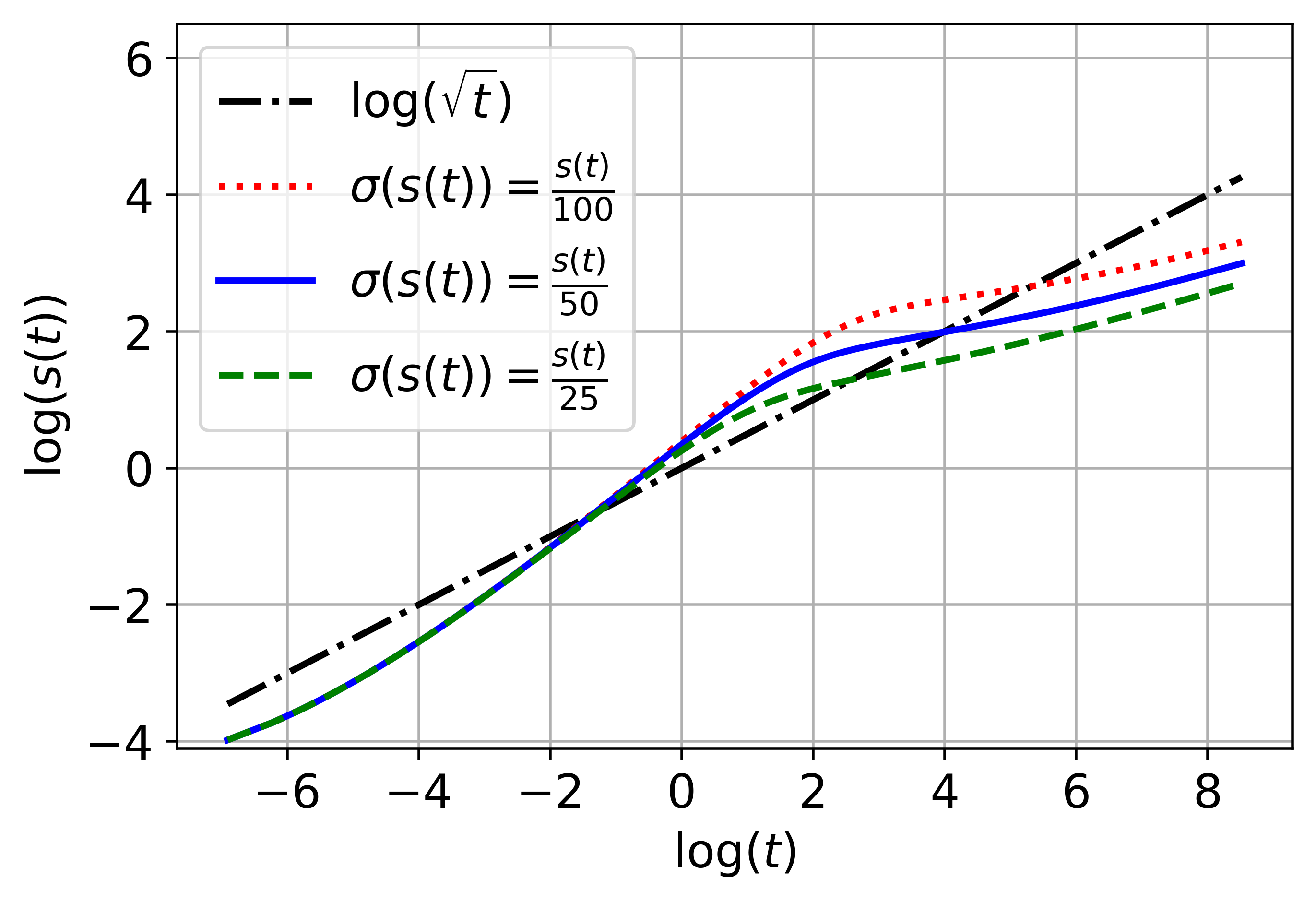}
	\hspace{0.1cm}
	\includegraphics[width=0.32\textwidth]{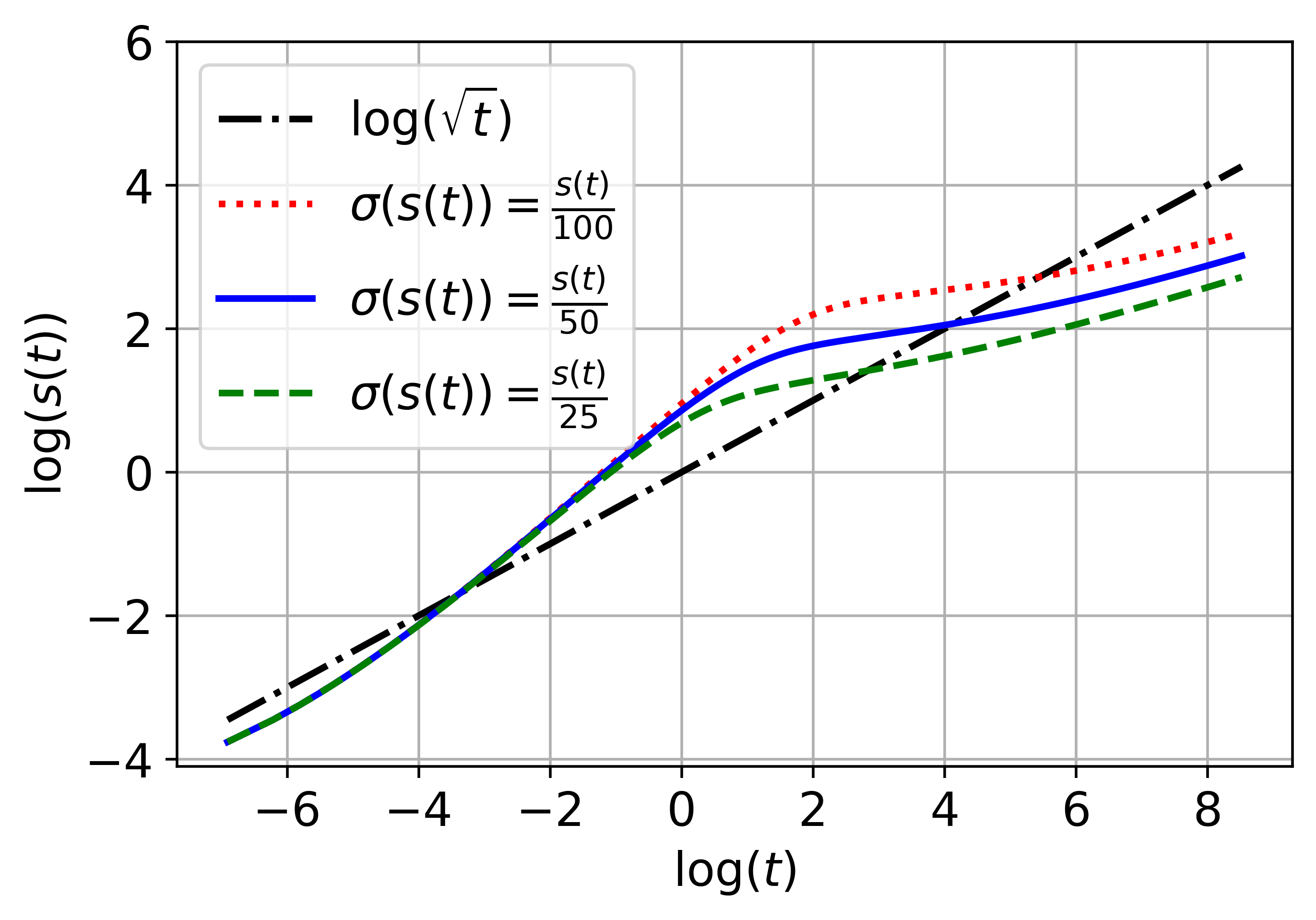}
	\caption{ For the foam rubber case: Comparison of $\log(\sqrt{t})$ \textit{vs}. $\log(t)$ with  $\log(s(t))$ \textit{vs}. $\log(t)$ for different values of $\sigma(s(t))$ with $T = 5000$ minutes and $a_0 = 500, \;a_0 =  1000, \;a_0 =  2000$ (from left to right).}
	\label{Fig:23}
\end{figure}

As we have already shown the concentration profile for large time in Figure \ref{Fig:4} (for the  dense rubber case) and Figure \ref{Fig:17} (for the foam rubber case), we now study the large time behaviour of diffusion fronts. We refer the reader to see our results for the dense rubber and foam rubber in Figure \ref{Fig:15} and Figure \ref{Fig:23}, respectively. We list in Table \ref{tab:2} and in Table \ref{tab:4} various best options for the exponent $\gamma$   matching the two choices of materials.

\begin {table}[h]
\begin{center}
		\begin{tabular}{ |p{2.5cm}|p{2cm}|p{2cm}| p{2cm}|}
			\hline
			\diagbox{$a_0$}{$\sigma(s(t))$} & $\frac{s(t)}{20}$&$\frac{s(t)}{10}$& $\frac{s(t)}{5}$\\
			\hline
			100  & 0.30145999& 0.2694511 & 0.23967876 \\
			500  &0.30531516 & 0.26961346 & 0.23571825\\
			1000&0.30841906& 0.27189137 &0.23634618 \\
			\hline
		\end{tabular}
		\caption {For the  dense rubber case: approximated value of  $\gamma$ entering $s(t) = t^\gamma$ for $T = 5000$ minutes.}
		\label{tab:2} 
	\end{center}
	\end {table}

\begin {table}[h]
\begin{center}
	\begin{tabular}{ |p{2.5cm}|p{2cm}|p{2cm}| p{2cm}|}
		\hline
		\diagbox{$a_0$}{$\sigma(s(t))$} & $\frac{s(t)}{100}$&$\frac{s(t)}{50}$& $\frac{s(t)}{25}$\\
		\hline
		500 & 0.39449518 & 0.35425142& 0.31698023\\
		1000  & 0.39890782& 0.35803523 & 0.32025865 \\
		2000  &0.40184547 & 0.3605881 & 0.32252626\\
		\hline
	\end{tabular}
	\caption {For the foam rubber case: approximated value of  $\gamma$ entering  $s(t) = t^\gamma$ for $T = 5000$ minutes.}
	\label{tab:4} 
\end{center}
\end {table}
Looking for instance at Figure \ref{Fig:6} (for the dense rubber case) or at Figure \ref{Fig:21} and Figure \ref{Fig:22} (for the foam rubber case), we notice that the penetration front $s(t)$ really behaves like $t^{\gamma}$ as $t$ is sufficiently large with $\gamma$ some positive number. We also notice that the dimensionless numbers Bi and $A_0$ are invariant with respect to time.  This makes us conjecture that, under suitable conditions,  there exist two positive constants $c_1$ and $c_2$ (independent of the choice of $t$)  such that the inequality
\begin{equation}\label{conjecture}
c_1t^{\gamma}\leq s(t)\leq c_2t^{\gamma}
\end{equation}
holds true for sufficiently large $t$. In Refs. \cite{aiki2011free,aiki2013large},  the authors were able to prove such estimate on the large time behavior of the moving boundary for a different set of  moving boundary conditions combined with a Dirichlet boundary condition at $x=0$. We expect that some of the techniques used in \cite{aiki2011free} 
are in principle applicable in our context as well. Hence,  we hope to be able to prove rigorously  the conjecture (\ref{conjecture}) in a follow-up paper. It is worth noting that if the conjecture were true, then new opportunities would appear particularly regarding the reduction of complex modeling approaches to the use of simple (material-dependent) educated guesses for fitted power laws to capture asymptotically the  penetration of the diffusants.

\section{Conclusion and outlook}\label{outlook}

Within this framework, we have proposed a one dimensional moving boundary model capable to capture the motion of sharp diffusion fronts in both dense and foam rubbers.  Including in our modeling work information about the measured swelling of the area and about the measured elongated height of the sample can be used further  to devise a balance law for the porosity change due to the penetration of the diffusant. This would implicitly include nonlinear diffusion effects. Such features need to be accounted for especially if one wants to capture swelling effects in more than one space dimension. Moreover, it is of both engineering and mathematical importance to build  simple asymptotic laws that can be verified experimentally in laboratory conditions and have a mathematically-controlled range of application outside the laboratory. From this point of view, it is important to understand the nature of the  deviations that usually occur in the position of the penetration front $s(t)$ from a $\sqrt{t}$-like behavior. As further possible outlet, we plan to  extend the current model to a multiscale framework, possibly based on ideas from \cite{aiki2020macromicro}. This approach would allow us  to include explicitly in our numerical investigations various microstructure models for the rubber foam.

\appendix

\section*{Notation}\label{A}


We list here the notation used in the text.\\
$A = 	N \times N$ tridiagonal matrix,\\
$A_0$ = Thiele modulus (dimensionless),\\
Bi = Biot number (dimensionless),\\
$C(0, T^*; V_k) =$ the space of continuous function from  $(0, T^*)$ to $V_k$  \cite{thomee1984galerkin},\\
$C^1(0,T^*) = $ the space of continuously differentiable function defined on $(0, T^*)$, \\
$D$ = diffusion constant for concentration in rubber (mm$^2$/min),\\
$H$ = Henry's constant (dimensionless),\\
$H^1(0,1)$ = Sobolev space \cite{thomee1984galerkin},\\
$K = N \times N$ tridiagonal matrix,\\ 
$M = N \times N$ positive definite matrix, \\
$V_k =$ finite dimensional subspace of $H^1(0,1)$ \cite{thomee1984galerkin},\\
$a_0$ = constant appearing in the speed of the moving front (mm$^4$/min/gram),\\
$b$ = concentration in the lower surface of the rubber (gram/mm$^3$),\\
$\textbf{e}_0,\; \textbf{e}_{N-1}$ = unit vectors,\\
$h(\tau)$ = position at time $\tau$ of the front separating a diffusant-free zone from a diffusant-penetrated zone within the rubber between concentration and rubber,  dimensionless form,\\
$m$= diffusant concentration (gram/mm$^3$),\\
$m_0$ = initial value of diffusant concentration $m$ (gram/mm$^3$),\\
$s(t)$ = position of the front between concentration and rubber at time $t$ (mm), \\
$s_0 = s(0)$ (mm), \\
$s^{\prime}(t)$ = growth rate of moving boundary (mm/min),\\
$u$ =  dimensionless form of diffusant concentration $m$,\\
$\beta$ = mass transfer constant from concentration to rubber at the boundary (mm/min),\\
$\phi_j$ = linear piece-wise continuous basis function.

\section*{Acknowledgements} A.M. thanks K. Kumazaki (Nagasaki, Japan) for very fruitful discussions on this topic. The work of S.N.  and A.M. is financed partly by the VR grant 2018-03648. T.A. and A.M. thank the Knowledge Foundation (project nr. KK 2019-0213) for supporting financially their research. The work of T.A. is partially supported also by JSPS KAKENHI project nr. JP19K03572.

\begin{center}
	\bibliographystyle{plain}
	\bibliography{polymerbib}
\end{center} 

\end{document}